\documentclass{elsart1p}
\usepackage{ifthen}
\usepackage{amsmath, latexsym, amssymb, graphics}
\usepackage[all, knot]{xy}
\input{grcawol.sty}

\newcommand{\gwnot}[2]{
   \grcalca = \grcolumn
   \multiply \grcalca by \factor
   \grcalcc = #1
   \multiply \grcalcc by \hfactor
   \advance \grcalca by \grcalcc
   \grcalcb = \grrow
   \multiply \grcalcb by \factor
   \advance \grcalcb by -\hfactor
   \put(\grcalca,\grcalcb) {\makebox(0,0){$\scriptstyle #2$}} }
\newcommand\gdnot[1]{\gwnot2{#1}}
\def\a{{\alpha}}
\def\b{{\beta}}
\def\e{{\epsilon}}
\def\w{{\omega}}
\date{}

\def\unit{{\mathbf{1}}}
\def\BQ{{\mathbb{Q}}}
\def\BC{{\mathbb{C}}}
\def\BF{{\mathbb{F}}}
\def\BZ{{\mathbb{Z}}}

\def\pf{\sprung{\it Proof.\,\,\,}}
\newcommand\triv{^{\operatorname{triv}}}
\newcommand\res{\operatorname{res}}

\def\Hom{{\mbox{\rm Hom}}}
\def\End{{\mbox{\rm End}}}
\def\Gal{{\mbox{\rm Gal}}}
\def\FSexp{{\mbox{\rm FSexp}}}

\newcommand\C[1]{{#1\mbox{-\bf{mod}}_{\operatorname{\mathsf fin}}}}

\newcommand\lcm{\operatorname{lcm}}

\newcounter{qedgemacht}
\newenvironment{proof}{\begin{pf}\setcounter{qedgemacht}{0}}{
\ifthenelse{\value{qedgemacht}=0}{\qed}{}\end{pf}}
\def\qedhere{\qed\stepcounter{qedgemacht}}

\numberwithin{thm}{section}
\numberwithin{equation}{section}

\newcommand\hh[1]{\widehat{#1}}
\newcommand\ord{\operatorname{ord}}

\newcommand\Z{\BZ}
\newcommand\ol[1]{\overline{#1}}
\newcommand\replace[1]{}
\newcommand\ptr{\operatorname{\underline{ptr}}}
\newcommand\pTr{\operatorname{\underline{ptr}}}
\newcommand\ptrl{\ptr^\ell}
\newcommand\ptrr{\ptr^r}

\newcommand\sprung{\smallskip}

\newcommand\sym{{\operatorname{sym}}}
\newcommand\id{\operatorname{id}}
\renewcommand\o{\otimes}

\newcommand\Tr{\operatorname{Tr}}

\newcommand\inv{^{-1}}
\DeclareMathOperator\ev{{\operatorname{ev}}}
\DeclareMathOperator\db{\operatorname{db}}
\newcommand\CC{\mathcal C}
\newcommand\CCstr{{\mathcal{C}_{\rm str}}}
\newcommand\DD{\mathcal D}
\newcommand\FF{\mathcal F}

\newcommand\du{^{\vee}}
\newcommand\bidu{^{\vee\vee}}
\newcommand{\leer}{\operatorname{--}}

\newcommand{\ou}[1]{\mathrel{\mathop{\otimes}_{#1}}}

\newcommand\catr{\underline{\operatorname{tr}}}
\newcommand\sw[1]{{}_{(#1)}}

\newcommand\so[1]{^{(#1)}}
\newcommand\som[1]{^{(-#1)}}
\newcommand\op{{\operatorname{op}}}
\newcommand\cop{{\operatorname{cop}}}

\makeatletter
\def\namelabel#1#2{\@bsphack
  \protected@write\@auxout{}%
         {\string\newlabel{#1.nme}{{#2}{#2}}}%
  \@esphack}
\def\nmlabel#1#2{\label{#2}\namelabel{#2}{#1}}
\newcommand\nmref[1]{\ref{#1.nme}\ \ref{#1}}
\makeatother

\begin{document}

\begin{frontmatter}
\title{Frobenius-Schur Indicators and Exponents of Spherical Categories}
\author{Siu-Hung Ng\thanksref{RichardsDank}}

\address{Department of Mathematics, Iowa State University, Ames, IA 50011, USA}

\thanks[RichardsDank]{The first author is supported by the NSA grant number
H98230-05-1-0020.}

\ead{rng@iastate.edu}

\thanks[PetersDank]{The second author is supported by the
\emph{Deutsche Forschungsgemeinschaft} through a Heisenberg
fellowship.}

\author{Peter Schauenburg\thanksref{PetersDank}}

\address{Mathematisches Institut der Universit\"at M\"unchen,
Theresienstr.\ 39, 80333 M\"unchen, Germany}

\ead{schauenburg@math.lmu.de}

\begin{abstract}
We obtain two formulae for the higher Frobenius-Schur indicators:
one for a spherical fusion category in terms of the twist of its
center and the other one for a modular tensor category in terms of
its twist. The first one is a categorical generalization of an
analogous result by Kashina, Sommerh\"{a}user, and Zhu for  Hopf
algebras, and the second one extends Bantay's 2nd indicator formula
for a conformal field theory to higher degree. These formulae imply
the sequence of higher indicators of an object in these categories
is periodic. We define the notion of Frobenius-Schur (FS-)exponent
of a pivotal category to be the global period of all these sequences
of higher indicators, and we prove that the FS-exponent of a
spherical fusion category is equal to the order of the twist of its
center. Consequently, the FS-exponent of a spherical fusion category
is a multiple of its exponent, in the sense of Etingof, by a factor
not greater than 2. As applications of these results, {we prove that
the exponent and the dimension of a semisimple quasi-Hopf algebra
$H$ have the same prime divisors, which answers two questions of
Etingof and Gelaki affirmatively for quasi-Hopf algebras. Moreover,
we prove that the FS-exponent of $H$ divides $\dim(H)^4$.} In
addition, if $H$ is a group-theoretic quasi-Hopf algebra, the
FS-exponent of $H$ divides $\dim(H)^2$, and this upper bound is
shown to be tight.
\end{abstract}
\end{frontmatter}
\section{Introduction}
We continue our investigation, begun in \cite{NS05,NS052}, of the
higher Frobenius-Schur indicators for quasi-Hopf algebras and, more
generally, certain fusion categories.\sprung

The classical (degree two) Frobenius-Schur indicator introduced a
century ago as well as the higher indicators are well-known
invariants of an irreducible representation of a finite group. The
degree two indicators for simple modules over a semisimple Hopf
algebra were studied by Linchenko and Montgomery \cite{LM00}, a
version for certain fusion categories by Fuchs, Ganchev, Szlachanyi,
and Vescernyes \cite{FGSV99}, and a more general version for simple
objects in pivotal categories by Fuchs and Schweigert
\cite{FucSch:CTCBC}. Bantay introduced a notion of Frobenius-Schur
indicator for a primary field of a conformal field theory via a
formula in terms of the modular data. The higher indicators for Hopf
algebras were introduced and studied in depth by Kashina,
Sommerh\"auser, and Zhu \cite{KSZ}. The degree two indicators for
simple modules of a semisimple quasi-Hopf algebra were studied by
Mason and Ng \cite{MN05}, and given a different treatment by
Schauenburg \cite{Sch04}. The higher indicators introduced in
\cite{NS05} are a generalization of all of the above to the case of
pivotal fusion categories, with more details and examples for the
special case of modules over a semisimple quasi-Hopf algebra worked
out in \cite{NS052}.\sprung

The definition of the indicator $\nu_n(V)$ of a simple object $V$ in
a pivotal fusion category $\CC$ in \cite{NS05} says that $\nu_n(V)$
is the trace of a certain endomorphism $E_V^{(n)}$ of $\CC(I,V^{\o
n})$. If $\CC=\C H$ for a semisimple Hopf algebra over $\BC$, we can
identify $\CC(I,V^{\o n})$ with the invariant subspace of $V^{\o
n}$, and $E_V^{(n)}$ corresponds to a cyclic permutation of tensor
factors. Thus, our definition corresponds to the ``first formula''
for the indicator given in \cite[Corollary 2.3]{KSZ}. The
\emph{definition} of the higher indicators in \cite{KSZ} is in terms
of Hopf powers of the integral, namely $\nu_n(V)=\chi(\Lambda\sw
1\cdots\Lambda\sw n)$, where $\Lambda$ is the normalized integral
and $\chi$ the character of $V$. This formula is generalized by the
description in \cite{NS05} of the $n$-th indicator as the pivotal
trace of the $n$-th Frobenius-Schur endomorphism of $V$; that this
is indeed a direct generalization of the defining formula in
\cite{KSZ} becomes evident from the calculations done in
\cite{NS052} for the quasi-Hopf algebra case.\sprung

In the present paper, we prove a generalization of the ``third
formula'' \cite[6.4, Corollary]{KSZ} for the higher indicators,
which says that the $n$-th indicator of $V$ is the trace of the
Drinfeld element of the double $D(H)$ taken on the induced module
$D(H)\ou HV$. The appropriate generalization to the categorical
setting is the following: For a spherical fusion category $\CC$,
consider the two-sided adjoint $K\colon\CC\to Z(\CC)$ of the
forgetful functor, where $Z(\CC)$ is the (left) center of $\CC$.
Then the $n$-th indicator $\nu_n(V)$ is $(\dim\CC)\inv$ times the
pivotal trace on $K(V)$ of the $n$-th power of the twist in the
ribbon category $Z(\CC)$. We prove this in section
\ref{sec:FS-ribbon}, making use of results of M\"uger on the adjoint
$K$ and Ocneanu's tube algebra for $\CC$. The key ingredient is
another formula for the $n$-th indicator in terms of the $n$-th
power of a special central element $t$ in the tube algebra of $\CC$
obtained in section \ref{sec:tube}.\sprung

In the Hopf algebra case, the ``second formula'' \cite[3.2,
Proposition]{KSZ}, of which the ``third formula'' is a reformulation
using the terminology of the Drinfeld double, connects the theory of
indicators to the exponent of a Hopf algebra studied by Kashina
\cite{Kas:OAHAHHYD,Kas:GPMHA} and Etingof and Gelaki \cite{EG99}. In
more detail, Kashina, Sommerh\"auser and Zhu define the exponent of
an irreducible representation of a semisimple complex Hopf algebra
$H$, and they characterize it as the order of a certain endomorphism
of the induced representation, the traces of whose powers are the
higher indicators. As a consequence of their results, the exponent
of $H$ is the least common multiple of the exponents of the
irreducible representations; the exponent of an irreducible
representation in turn is the period of the sequence formed by its
higher indicators, or the least number $n$ such that the $n$-th
indicator is the representation's dimension.\sprung

Theorem \ref{t:nu_tr_formula} implies that the sequence
$\{\nu_n(V)\}_n$ of the higher indicators  of an object $V$ in the
spherical fusion category $\CC$ is periodic as well. Moreover, the
$m$-th term of the sequence is the pivotal dimension $d(V)$ of $V$
whenever $m$ is a multiple of the order of the twist of
$Z(\CC)$.\sprung

    In section \ref{sec:FSexp} we define the Frobenius-Schur exponent of
an object $V$ in a pivotal category $\CC$ to be the least positive
integer $n$ such that $\nu_n(V)$ is the pivotal dimension $d(V)$ of
$V$, and we define the Frobenius-Schur exponent of $\CC$ to be the
least positive integer $n$ such that $\nu_n(V)=d(V)$ for all $V \in
\CC$. We then prove that if $\CC$ is a spherical fusion category,
then its Frobenius-Schur exponent is equal to the order of the twist
of $Z(\CC)$. By \cite{Etingof02} the exponent is finite and our
result implies that the exponent of $\CC$ divides the
Frobenius-Schur exponent. As it turns out in the first examples of
the same section, the Frobenius-Schur exponent is different in
general from the exponent as defined by Etingof \cite{Etingof02}. In
section \ref{sec:versus}, however, we show that the Frobenius-Schur
exponent can at most be twice  the exponent. \sprung

Bantay introduced the (degree 2) Frobenius-Schur indicator for a
primary field $V$ of a conformal field theory via a formula in terms
of the modular data of the CFT and he showed that the indicator of
$V$ is 0 if $V \not\cong V^*$ and $\pm 1$ if $V^* \cong V$ (cf.
\cite{Bantay97}). In section \ref{s:MTC}, we derive a formula for
the $n$-th indicator of a simple object $V$ of a modular tensor
category by computing the trace of $E_V^{(n)}$. Our formula for
$\nu_n(V)$ contains Bantay's formula for degree two indicators as
the special case $n=2$. An important consequence of this formula is
the invariance of Frobenius-Schur exponent of a spherical fusion
category under the center construction. Moreover, the
Frobenius-Schur exponent of a modular tensor category is equal to
the order of the twist.\sprung

In section \ref{sec:Cauchy} we obtain several results for complex
semisimple quasi-Hopf algebras $H$, which have been known for Hopf
algebras, making use of Frobenius-Schur indicators and exponents. We
prove that the dimension of $H$ is even if $H$ admits a self-dual
simple module, and we generalize the Hopf algebra version of
Cauchy's Theorem from \cite{KSZ} to our setting: The exponent, {in
the sense of Etingof \cite{Etingof02},} and the dimension of a
quasi-Hopf algebra $H$ have the same prime factors. {This result
answer two questions of Etingof and Gelaki \cite{EG99} affirmatively
for semisimple complex quasi-Hopf algebras. We have also shown that
the exponent and the Frobenius-Schur exponent coincide if the
dimension of $H$ is odd.}  \sprung

In section \ref{sec:bounds} we prove two bounds for the
Frobenius-Schur exponent. In \cite{Kas:GPMHA} Kashina asked whether
the exponent of a semisimple complex Hopf algebra always divides the
dimension. Etingof and Gelaki \cite[Theorem 4.3]{EG99} have shown
that the exponent divides the third power of the dimension. As an
immediate consequence of the results in \cite{Etingof02} the
Frobenius-Schur exponent of a semisimple complex quasi-Hopf algebra
$H$ divides the fifth power of the dimension of $H$. We improve this
bound by lowering the exponent to the fourth power. We also study an
important class of quasi-Hopf algebras, namely the group-theoretical
quasi-Hopf algebras corresponding to the group-theoretical fusion
categories introduced by Ostrik \cite{Ostrik03}. For a
group-theoretical quasi-Hopf algebra $H$ we can show that the
(Frobenius-Schur) exponent divides the square of the dimension. More
precisely, we can express the Frobenius-Schur exponent of a
group-theoretical quasi-Hopf algebra, which is constructed from a
finite group $G$ and a three-cocycle $\omega$ on $G$, in terms of
the cohomology class of $\omega$ and its restrictions to the cyclic
subgroups of $G$. This same description of the exponent and the
resulting general bound on the exponent for group-theoretical
categories were recently obtained by Natale \cite{Nat05}.

\section{Preliminaries} \label{s:prelim}
We will collect some conventions and facts on monoidal categories
and quasi-Hopf algebras. Most of these are well-known (and we refer
to \cite{NS05,NS052} and the literature cited there), the others are
easy observations.\sprung

In a monoidal category $\CC$ with tensor product $\o$ we denote the
associativity isomorphism by $\Phi\colon (U\o V)\o W\to U\o(V\o W)$.
We assume that the unit object $I\in\CC$ is strict. If $X,Y\in\CC$
are obtained by tensoring together the same sequence of objects with
two different arrangements of parentheses, there is an isomorphism
between them which is obtained by composing several instances of
$\Phi$ or $\Phi\inv$; it is unique by coherence, and will be denoted
by $\Phi^?\colon X\to Y$. A monoidal functor $\mathcal
F\colon\CC\to\mathcal D$ preserves tensor products by way of a
coherent isomorphism $\xi\colon\mathcal F(V)\o\mathcal
F(W)\to\mathcal F(V\o W)$ and $\mathcal F(I)=I$. An equivalence of
monoidal categories is a monoidal functor that is an
equivalence.\sprung

A left dual object $(V\du,\ev,\db)$ consists of an object $V\du$ and
morphisms $\ev\colon V\du\o V\to I$ and $\db\colon I\to V\o V\du$
such that
\begin{gather*}
V\xrightarrow{\db\o V}(V\o V\du)\o V\xrightarrow{\Phi}V\o(V\du\o
V)\xrightarrow{V\o\ev}V,
\\
V\du\xrightarrow{V\du\o\db}V\du\o(V\o
V\du)\xrightarrow{\Phi\inv}(V\du\o V)\o V\du\xrightarrow{\ev\o
V\du}V\du
\end{gather*}
are identities. Right duals are defined analogously. If every object
has a (left) dual, $\CC$ is called (left) rigid. If $\CC$ is left
rigid, taking duals extends naturally to a monoidal functor
$(\leer)\du\colon\CC^\op\to\CC$. Double dualization is consequently
a monoidal functor $(\leer)\bidu\colon\CC\to\CC$. A pivotal monoidal
category is a left rigid monoidal category equipped with an
isomorphism $j\colon V\to V\bidu$ of monoidal functors. Let $f\colon
V\to V$ be a morphism in a pivotal category $\CC$. The left and
right pivotal traces of $f$ are
\begin{gather*}
\ptrr(f)=\catr(j_Vf)=\left(I\xrightarrow{\db}V\o
V\du\xrightarrow{f\o V\du}V\o V\du\xrightarrow{j_V\o V\du}V\bidu\o
V\du\xrightarrow{\ev}I\right),
\\
\ptrl(f)=\left(I\xrightarrow{\db}V\du\o V\bidu\xrightarrow{V\du\o
j_V\inv}V\du\o V\xrightarrow{V\du\o f}V\du\o
V\xrightarrow{\ev}I\right).
\end{gather*}
The left and right pivotal dimensions of $V\in\CC$ are
$d_\ell(V)=\ptr^\ell(\id_V)$ and $d_r(V)=\ptr^r(\id_V)$. If the left
and right traces of every morphism are the same, then $\CC$ is
called a spherical monoidal category. In this case, traces and
pivotal dimensions will be denoted by $\ptr(f)$ and $d(V)$. If $\CC$
is $\BC$-linear, semisimple, and $d_\ell(V)=d_r(V)$ for each simple
$V$, then $\CC$ is spherical.\sprung

Any monoidal category is equivalent as a monoidal category to a
strict monoidal category, that is, one in which the associativity
isomorphism $\Phi$ is the identity. A pivotal monoidal category
$\CC$ is, moreover, equivalent as a pivotal monoidal category to a
strict pivotal monoidal category $\CCstr$, that is, a pivotal
monoidal category in which the associativity isomorphism, the
pivotal structure $j$, and the canonical isomorphism $(V\o W)\du \to
W\du \o V\du$ are identities. Equivalence as pivotal monoidal
categories means that the monoidal equivalence $\CC\to \CCstr$
preserves pivotal structures in a suitable sense; we refer to
\cite{NS05} for details. If $\CC$ is spherical, then so is
$\CCstr$.\sprung

In a strict monoidal category we make free use of graphical
calculus. For instance, the condition on a pivotal category to be
spherical is depicted as
$$
  \ptrr(f)=
  \gbeg35
  \gwdb{3} \gnl
  \gbmp{f}\gvac1\gcl1 \gnl
  \gcl1 \gvac1 \gcl1 \gnl
  \gbmp{j}\gvac1\gcl1 \gnl
  \gwev{3}\gnl
  \gend
  =
  \gbeg35
  \gwdb{3} \gnl
  \gcl1\gvac1 \gbmp{\,j^{-\!\!1}} \gnl
  \gcl1 \gvac1 \gcl1 \gnl
  \gcl1\gvac1 \gbmp{f} \gnl
  \gwev{3}\gnl
  \gend = \ptrl(f)\,.
$$
The strict pivotal case allows us to simply drop any instance of $j$
resulting in
$$ \ptr(f) =
  \gbeg23
  \gdb \gnl
  \gbmp{f}\gcl1 \gnl
  \gev\gnl
  \gend
  =
  \gbeg23
  \gdb \gnl
  \gcl1\gbmp{f} \gnl
  \gev\gnl
  \gend\,.
$$

Now let $\CC$ be a braided monoidal category with braiding $c$. The
Drinfeld isomorphism is the natural isomorphism $u: Id \rightarrow
(-)\bidu$ given by
$$
  u_V:= (\ev_V\o V\bidu)\circ \Phi\inv_{V\du, V, V\bidu} \circ
  (V\du \o c_{V,V\bidu}\inv) \circ
  \Phi_{V\du, V\bidu, V}\circ
  (\db_{V\du}\o V).
$$
If $\CC$ is strict, then we have the graphical representation
$$
  u_V=\gbeg35
  \gvac2\got1V\gnl
  \gdb\gcl1 \gnl
  \gcl1\gibr\gnl
  \gev\gcl1\gnl
  \gvac2 \gob1{V\bidu}
  \gend,
  \qquad\text{ where }
  c_{VW}=\gbeg23\got1V\got1W\gnl\gbr\gnl\gob1W\gob1V\gend
  \text{ and }
  (c_{VW})\inv=\gbeg23\got1W\got1V\gnl\gibr\gnl\gob1V\gob1W\gend\,.
$$
We see that $u_{V\o W}=(u_V\o u_W)c_{VW}\inv c_{WV}\inv$. In
particular, there is a bijective correspondence between pivotal
structures $j$ on $\CC$ and {\it twists}, that is, automorphisms
$\theta$ of the identity endofunctor satisfying
\begin{equation} \label{def:twist}
\theta_{V\o W}=(\theta_V\o\theta_W)c_{WV}c_{VW} \quad
\text{and}\quad \theta_I=\id_I,
\end{equation}
given by $\theta=u\inv j$.\sprung

The (left) center $Z(\CC)$ of a monoidal category $\CC$ has objects
pairs $(V, e_V)$ with $V \in \CC$ and $e_{V}(-): V \o (-)\;
\rightarrow \; (-) \o V$ a natural isomorphism satisfying the
properties $e_V(I)=\id_V$ and
$$
(X\o e_V(Y))\circ \Phi_{X,V,Y} \circ (e_V(X) \o Y)  = \Phi_{X,Y,
V}\circ e_V(X \o Y)\circ \Phi_{V,X,Y}
$$
for all $X, Y \in \CC$. The tensor product of $Z(\CC)$ is given by $
  (V, e_V)\o (W, e_W) = (V \o W, e_{V \o W}),
$ where
$$
e_{V \o W}(X) = \Phi_{X, V, W} \circ (e_V(X) \o W)\circ  \Phi_{V, X,
W}\inv \circ (V \o e_W(X)) \circ \Phi_{V, W, X}
$$
for any $X \in \CC$, and the neutral object is $(I, e_I)$ with
$e_I(X)=\id_X$. The associativity isomorphism in $Z(\CC)$ is that of
$\CC$. In this way $Z(\CC)$ is a monoidal category, and braided with
braiding given by $e_V(W)$.\sprung

If $\CC$ is left rigid, then $Z(\CC)$ is rigid; the dual of $(V,
e_V) \in Z(\CC)$, is $(V\du,e_{V\du})$ with
$$
e_{V\du}(X)= ( \ev_V \o (X \o V\du))\circ \Phi^? \circ (V\du \o
(e_V(X)\inv \o V))\circ \Phi^?\circ (V\du \o (X \o \db_V))\,.
$$
The evaluation and dual basis  homomorphisms for the dual in
$Z(\CC)$ are those of the dual $V\du$ in $\CC$. If $\CC$ is pivotal,
the pivotal structure $j: Id \rightarrow (-)\bidu$ induces a pivotal
structure in $Z(\CC)$. If $\CC$ is spherical, so is $Z(\CC)$.\sprung

Any equivalence $\FF\colon\CC\to\mathcal D$ of monoidal categories
induces in a natural way an equivalence $\hat\FF\colon Z(\CC)\to
Z(\mathcal D)$ of braided monoidal categories. In addition, if $\CC$
and $\DD$ are pivotal monoidal categories and $\FF$ preserves their
pivotal structures, then $\hat\FF\colon Z(\CC)\to Z(\mathcal D)$
also preserves their pivotal structures. Moreover, $\hat \FF$
preserves the twists $\theta$ associated with their pivotal
structures, i.e.
$$
\hat{\FF}(\theta_{(V,e_V)}) = \theta_{\hat{\FF}(V, e_V)}\,.\sprung
$$

A fusion category over $\BC$ is a rigid $\BC$-linear monoidal
category which is semisimple with finitely many nonisomorphic simple
objects, whose endomorphism rings are isomorphic to $\BC$. If $\CC$
is a strict spherical fusion category with a braiding $c$, then by
\eqref{def:twist} the twist $\theta$ associated with the pivotal
structure of $\CC$ is identical to $u\inv$ where $u$ is the Drinfeld
isomorphism associated with the braiding of $Z(\CC)$. Moreover, for
any simple object $V$ of $\CC$,
$$
\ptr\left(\theta_{V}\right) = \ptr\left( \gbeg35
  \got1{\,\,V}\gvac2\gnl
  \gcl1\gdb \gnl
  \gbr\gcl1 \gnl
  \gcl1\gev\gnl
  \gob1{\,\,V}\gvac2
  \gend
\right) = \gbeg43
  \gdb \gdb \gnl
  \gcl1 \gbr\gcl1 \gnl
  \gev\gev\gnl
  \gend =
  \ptr\left(
\gbeg35
  \gvac2\got1{V\,}\gnl
  \gdb\gcl1 \gnl
  \gcl1\gbr \gnl
  \gev\gcl1\gnl
  \gvac2 \gob1{V}
  \gend
\right)
$$
and so
$$
 \theta_{V} = u_V\inv
 =\gbeg35
  \gvac2\got1V\gnl
  \gdb\gcl1 \gnl
  \gcl1\gbr \gnl
  \gev\gcl1\gnl
  \gvac2 \gob1{V}
  \gend
  =\gbeg35
  \got1V\gvac2\gnl
  \gcl1\gdb \gnl
  \gbr\gcl1 \gnl
  \gcl1\gev\gnl
  \gob1{V}\gvac2
  \gend, \quad\text{and}\quad
  \theta_{V} \inv
 =\gbeg35
  \gvac2\got1V\gnl
  \gdb\gcl1 \gnl
  \gcl1\gibr \gnl
  \gev\gcl1\gnl
  \gvac2 \gob1{V}
  \gend
  =\gbeg35
  \got1V\gvac2\gnl
  \gcl1\gdb \gnl
  \gibr\gcl1 \gnl
  \gcl1\gev\gnl
  \gob1{V}\gvac2
  \gend\,.
$$
In particular, $\CC$ is a ribbon category with respect to the twist
$\theta$.\sprung

If $H$ is a quasi-Hopf algebra with associator $\phi\in H^{\o 3}$
and (quasi-)antipode $(S,\alpha,\beta)$, then the category $\C H$ of
finite-dimensional left $H$-modules is a rigid monoidal category
with associativity isomorphism $\Phi$ given as left multiplication
with $\phi$, dual object $V\du=V^*$ the vector space dual with
module structure the transpose of the action through $S$, evaluation
$\ev(f\o v)=f(\alpha v)$, and dual basis morphism $\db(1)=\sum_i
\beta v_i\o v^i$, where $\{v_i\}_i$ and $\{v^i\}_i$ are dual bases
for $V$ and $V^*$ respectively. If $H$ is semisimple, we denote the
normalized integral in $H$ by $\Lambda$.\sprung

By a result of Etingof, Nikshych and Ostrik \cite[Section 8]{ENO},
the categories $\C H$ for semisimple quasi-Hopf algebras $H$ can be
characterized as those fusion categories $\CC$ for which the
Frobenius-Perron dimension of every simple object is an integer. For
such a category \cite{ENO} show that $\CC$ is pseudo-unitary, and it
has a pivotal structure $j$ determined by the condition $d(V)=\ev_V
(j_V\o \id) \db_V=\dim(V)\id_\BC$. We will call this the canonical
pivotal structure; it is spherical. For $\CC=\C H$ with a quasi-Hopf
algebra $H$ we have $j(v)=g\inv v$ for a certain element $g$ of $H$
called the trace-element of $H$. \sprung

\begin{rem}\label{r:piv_dim}
 If $j'$ is a second pivotal structure on a pivotal fusion category
$\CC$ over $\BC$ with pivotal structure $j$, we have $j'=j\lambda$
for some monoidal automorphism $\lambda$ of the identity functor.
For a simple $V$, the component $\lambda_V$ is a scalar. In
particular, $\lambda_I=1$. If we denote the left and right pivotal
dimensions with respect to the pivotal structure $j'$ by $d'_\ell$
and $d'_r$, we see $d_\ell'(V)=\lambda_V\inv d_\ell(V)$ and
$d_r'(V)=\lambda_V d_r(V)$ for simple $V$. Also, since $V\o V\du$
contains $I$, we have $\lambda_{V\du}\lambda_V=1$. By
\cite[Proposition 2.9]{ENO},
$$
|d_\ell'(V)|^2 =|V|^2 = |d_\ell(V)|^2\,.
$$
In particular, the absolute values of the pivotal dimensions of
simple objects and the pivotal dimension of $\CC$ do not depend on
the choice of a pivotal structure. Finally, if both pivotal
structures are spherical, then $\lambda_V=\pm 1$ for each simple
$V$. \sprung
\end{rem}

\begin{rem}\label{rem:rightcenter}
  Analogous to the left center construction,
 the right version of the center $\ol{Z}(\CC)$ of $\CC$
consists of objects pairs $(V, \bar{e}_V)$ with $V \in \CC$ and
$\bar{e}_{V}(-): (-) \o V \rightarrow V \o (-)$ a natural
isomorphism satisfying the properties $\bar{e}_V(I)=\id_V$ and
$$
\bar{e}_V(X \o Y)= \Phi_{X,Y, V}\circ(\bar{e}_V(X)\o Y)\circ
\Phi\inv_{X,V,Y} \circ (X \o \bar{e}_V(Y) ) \circ \Phi_{X,Y,V}
\quad\text{for } X, Y \in \CC\,.
$$
The tensor product of $\ol{Z}(\CC)$ is given by $
  (V, \bar{e}_V)\o (W, \bar{e}_W) = (V \o W, \bar{e}_{V \o W}),
$ where
$$
\bar{e}_{V \o W}(X) =   \Phi_{V, W, X}\inv\circ (V \o
\bar{e}_W(X))\circ  \Phi_{V, X, W} \circ (\bar{e}_V(X) \o W) \circ
\Phi_{X, V, W}\inv
$$
for any $X \in \CC$, and the neutral object is $(I, \bar{e}_I)$ with
$\bar{e}_I(X)=\id_X$. The associativity isomorphism in $\ol{Z}(\CC)$
is also the same as in $\CC$. The right center $\ol{Z}(\CC)$ is also
a braided monoidal category with braiding given by $\bar{e}_V(W):
(W, \bar{e}_W) \o (V, \bar{e}_V) \to  (V, \bar{e}_V)\o (W,
\bar{e}_W)$. In addition, if $\CC$ is left rigid, then $\ol{Z}(\CC)$
is rigid; the dual of $(V, \bar{e}_V) \in \ol{Z}(\CC)$, is
$(V\du,\bar{e}_{V\du})$ with
$$
\bar{e}_{V\du}(X)\inv = (\ev_V \o (X \o V\du))\circ \Phi^? \circ
(V\du \o (\bar{e}_V(X) \o V))\circ \Phi^?\circ (V\du \o (X \o
\db_V))
$$
$$
\ev_{(V, \bar{e}_V)} = \ev_V, \quad \text{and} \quad \db_{(V,
\bar{e}_V)} = \db_V\,.
$$

By \cite{JS93}, the natural isomorphism $c'_{(V, e_V),(W, e_W)} :=
e_W(V)\inv$ for $(V, e_V)$, $(W, e_W) \in Z(\CC)$ also defines a
braiding on $Z(\CC)$, and we denote by $Z'(\CC)$ the braided
monoidal category $Z(\CC)$ with the braiding $c'$. Then $Z'(\CC)$
and $\ol{Z}(\CC)$ are equivalent braided monoidal categories under
monoidal equivalence $(T, \xi): Z'(\CC) \to \ol{Z}(\CC)$ with $\xi:
T(V, e_V) \o T(W, e_W) \to T(V\o W, e_{V \o W})$  the identity,
$T(V, e_V) = (V, e_V\inv)$ and $T(f)=f$ for any objects $(V, e_V)$,
$(W, e_W)$ and map $f$ of $Z'(\CC)$.  If $\CC$ admits a pivotal
structure, $(T, \xi)$ preserves the induced pivotal structures of
$Z'(\CC)$ and $\ol Z(\CC)$ as well as their associated twists.
\sprung

In addition, if $\CC$ is a spherical fusion category over $\BC$ with
the pivotal structure $j$, then $Z(\CC)$, $Z'(\CC)$ and
$\ol{Z}(\CC)$ are pivotal fusion categories with their pivotal
structures inherited from $\CC$. Let $\theta, \theta'$ and
$\bar{\theta}$ be the associated twists of $Z(\CC)$, $Z'(\CC)$ and
$\ol{Z}(\CC)$ respectively and $(V,e_V)$  a simple object of
$Z(\CC)$. Then $\theta_{(V, e_V)}= \w \id_{(V, e_V)}$ for some
non-zero scalar $\w$ in $\BC$. By considering the strictifications
of these pivotal fusion categories, we have $\theta'_{(V,e_V)} =
\w\inv\id_{(V, e_V)}$ and hence $\bar{\theta}_{T(V,e_V)}=\w\inv
\id_{T(V, e_V)}$. In particular, we have
$$
\theta'_{(V, e_V)} = \theta\inv_{(V, e_V)}\,.
$$
\end{rem}

\section{Tube Algebra of a Strict Spherical Fusion
Category}\label{sec:tube} In this section, we consider Ocneanu's
tube algebra $\Theta_L$ of a strict spherical fusion category and a
special element $t \in \Theta_L$. This element $t$ has been
considered in \cite{Izumi00} and \cite{MugerII03} for the
computations of the Gauss sums of $\CC$. We will show in Lemma
\ref{l:phi_t^n} that the $n$-th FS-indicator can be expressed in
terms of $t^n$. This observation is essential to our proofs of
Theorem \ref{t:nu_tr_formula} and Proposition \ref{t:dD_relation}.
\sprung

Let $\CC$ be a strict spherical fusion category over $\BC$. Let
$X_i$, $i \in \Gamma$, the set of isomorphism classes of simple
objects of $\CC$. Since the left dual and right dual of an object in
$\CC$ are isomorphic, we simply denote the left(or right) dual of
any object $V \in \CC$ by $V^*$. For any $i \in \Gamma$, we define
$\bar{i}$ by the equation
$$
X_i^* = X_{\bar{i}}
$$
and we define $0\in\Gamma$ by $X_0=I$, the neutral object of $\CC$.
\sprung

For any $V \in \CC$, we let
$$
d(V) = \ptr(\id_V), \quad d_i = d(X_i)\quad \mbox{for }i \in \Gamma,
\quad \mbox{and}\quad\dim \CC = \sum_{i \in \Gamma} d_i d_{\bar{i}}.
$$
Note that $d_0=1$ and $d_i = d_{\bar{i}}$ for $i \in \Gamma$.
Following \cite{MugerII03}, we consider Ocneanu's tube algebra
$\Theta_L=\bigoplus_{i, j, k \in \Gamma} \CC(X_i \o X_j, X_j \o
X_k)$ with multiplication given by
$$
uv[i,j,k]=\frac{1}{d_j\lambda} \sum_{s, l,m \in \Gamma} d_l
d_m\sum_{\a=1}^{N_{lm}^j}\left(\, \gbeg{5}{11}
  \got1{X_i}\gvac2\got1{X_j} \gnl
  \gcl1\gvac2\gcl1\gnl
  \gcl1\gvac1\gdnot{\quad p^\a_{j, lm}}\glmpb\gcmpt\grmpb\gnl
  \gcl1\gvac1\gcl1\gvac1\gcl1\gnl
  \gdnot{\quad v[i,l,s]}\glmptb\gcmp\grmptb\gvac1\gcl1\gnl
  \gcl1\gvac1\gcl1\gvac1\gcl1\gnl
  \gcl1\gvac1\gdnot{\quad u[s,m,k]}\glmptb\gcmp\grmptb\gnl
  \gcl1\gvac1\gcl1\gvac1\gcl1\gnl
  \gdnot{\quad q^\a_{lm, j}}\glmpt\gcmpb\grmpt\gvac1\gcl1\gnl
  \gvac1\gcl1\gvac2\gcl1\gnl
  \gvac1\gob1{X_j}\gvac2\gob1{X_k}
  \gend\,
\right)
$$
where $\{p^\a_{j, lm}\}_\a$ is a basis for $\CC(X_j, X_l\o X_m)$ and
$\{q^\a_{lm,j}\}_\a$ the basis for $\CC(X_l\o X_m, X_j)$ dual to
$\{p^\a_{j, lm}\}_\a$ and $\lambda = \sqrt{\dim \CC}$. The product
of $\Theta_L$ is independent of the choice of basis for $\CC(X_j,
X_l\o X_m)$ (cf. \cite{Ocneanu93}\cite{Izumi00}\cite{EvKa95} for the
original definition of Ocneanu's tube algebra). The identity element
$\unit$ of $\Theta_L$ is given by
$$
\unit[i,j,k]=\lambda \delta_{j,0}\delta_{i,k}\id_{X_i}\,.
$$

Let $\Theta_i = \bigoplus_{j \in \Gamma} \CC(X_i \o X_j, X_j \o
X_i)$ and $\Theta_\CC=\bigoplus_{i \in \Gamma} \Theta_i$. For any $u
\in \Theta_\CC$, one can define $\hat{u} \in \Theta_L$ given by
$$
\hat{u}[i,j,k]= \delta_{i,k} u[i,j] \quad\mbox{for all } i,j,k \in
\Gamma.
$$
We will identify $\Theta_\CC$ with a subspace of $\Theta_L$ under
this identification. It is easy to see that $\Theta_\CC$ is closed
under the multiplication on $\Theta_L$, and it also contains the
identity element $\unit$ of $\Theta_L$. Consider the element $t \in
\Theta_\CC$ given by
\begin{equation}\label{eq:t}
 t[i,j] = \frac{\lambda}{d_i} \delta_{ij} \id_{X_i^{\o 2}}\,.
\end{equation}
\sprung
\begin{lem}\label{l:3.1}
  \begin{equation}\label{eq:t^n}
  t^n [i,j] =\frac{\lambda}{d_j} \sum_{\a=1}^{N_{i^n}^j}\left(\,
  \gbeg{3}{7}
  \got1{X_i}\gvac1\got1{X_j} \gnl
  \gcl1\gvac1\gcl1\gnl
  \gcl1\gdnot{f^\a_{j, i^n}}\glmpb\grmptb\gnl
  \gcl1\gcl1\gcl1\gnl
  \gdnot{g^\a_{i^n, j}}\glmptb\grmpt\gcl1\gnl
  \gcl1\gvac1\gcl1\gnl
  \gob1{X_j}\gvac1\gob1{X_i}
  \gend\,
\right)
  \end{equation}
  for $n \ge 1$ where $N_{i^n}^j=
  \dim \CC(X_j, X^{\o i})$, $\{f^\a_{j, i^n}\}_\a$ is a basis for
  $\CC(X_j, X_i^{\o n})$, and $\{g^\a_{i^n, j}\}_\a$ is the dual basis for $\CC(X_i^{\o n},
  X_j)$.
\end{lem}
\begin{proof}
For $n=1$, the formula holds by definition of $t$. Assume that $t^n$
is given by \eqref{eq:t^n}. Then
$$
t^{n+1}[i,j]=\frac{1}{d_j\lambda} \sum_{l,m \in \Gamma} d_l
d_m\sum_{\a=1}^{N_{lm}^j}\left(\, \gbeg{3}{11}
  \got1{X_i}\gvac1\got1{X_j} \gnl
  \gcl1\gvac1\gcl1\gnl
  \gcl1\gdnot{p^\a_{j, lm}}\glmpb\grmptb\gnl
  \gcl1\gcl1\gcl1\gnl
  \gdnot{t^n}\glmptb\grmptb\gcl1\gnl
  \gcl1\gcl1\gcl1\gnl
  \gcl1\gdnot{t}\glmptb\grmptb\gnl
  \gcl1\gcl1\gcl1\gnl
  \gdnot{q^\a_{lm, j}}\glmptb\grmpt\gcl1\gnl
  \gcl1\gvac1\gcl1\gnl
  \gob1{X_j}\gvac1\gob1{X_i}
  \gend\,
\right) =\frac{\lambda}{d_j} \sum_{l\in \Gamma}
\sum_{\b=1}^{N_{i^n}^j}\sum_{\a=1}^{N_{li}^j}\left(\,
  \gbeg{4}{11}
  \got1{X_i}\gvac2\got1{X_j} \gnl
  \gcl1\gvac2\gcl{1}\gnl
  \gcl1\gvac1\gdnot{p^\a_{j, li}}\glmpb\grmptb\gnl
  \gcl1\gvac1\gcl1\gcl1\gnl
   \gcl1\gdnot{f^\b_{l, i^n}}\glmpb\grmptb\gcl1\gnl
  \gcl1\gcl1\gcl1\gcl1\gnl
  \gdnot{g^\b_{i^n, l}}\glmptb\grmpt\gcl1\gcl1\gnl
  \gcl1\gvac1\gcl1\gcl1\gnl
  \glmptb\gdnot{q^\a_{li, j}}\gcmp\grmpt\gcl1\gnl
  \gcl1\gvac2\gcl1\gnl
  \gob1{X_j}\gvac2\gob1{X_i}
  \gend\,
\right)
$$
Note that  $\{(f^\b_{l,i^n} \o X_i)\circ p^\a_{j, li}\}_{\a, \b, l}$
forms a basis for $\CC(X_j, X_i^{\o (n+1)})$ with dual basis
$$
\{q^\a_{li, j}\circ (g^\b_{i^n, l} \o X_i)\}_{\a, \b, l}
$$
for $\CC(X_i^{\o (n+1)}, X_j)$. Thus, the result follows by
induction.
\end{proof}
Define $\phi$ and $\phi_i \in \Theta_\CC^*$ by
$$
\phi(u) = \lambda \sum_{i\in \Gamma} d_i \pTr_{X_i}(u[i,0])\,,\quad
\mbox{and}\quad  \phi_i(u) = \phi\circ \pi_i(u)
$$
where $\pi_i$ is the natural projection from $\Theta_\CC$ to
$\Theta_i$.\sprung

\begin{lem}\label{l:phi_t^n}
For any $i \in \Gamma$, $\lambda^2 d_i
\overline{\nu_n(X_i)}=\phi_i(t^n)$.
\end{lem}
\begin{proof}
By Lemma \ref{l:3.1},
$$
\begin{aligned}
  \phi_i(t^n)& =\frac{\lambda^2}{d_0} \sum_{\a=1}^{N_{i^n}^0}d_i\pTr_{X_i}\left(\,
  \gbeg{3}{7}
  \got1{X_i}\gvac1\got1{I} \gnl
  \gcl1\gvac1\gcl1\gnl
  \gcl1\gdnot{f^\a_{I, i^n}}\glmpb\grmptb\gnl
  \gcl1\gcl1\gcl1\gnl
  \gdnot{g^\a_{i^n, I}}\glmptb\grmpt\gcl1\gnl
  \gcl1\gvac1\gcl1\gnl
  \gob1{I}\gvac1\gob1{X_i}
  \gend\,
\right) =\lambda^2 \sum_{\a=1}^{N_{i^n}^0}d_i\left(\,
  \gbeg{4}{5}
  \gwdb{4} \gnl
  \gcl1\gdnot{f^\a_{0, i^n}}\glmpb\grmpb\gcl1\gnl
  \gcl1\gcl1\gcl1\gcl1\gnl
  \gdnot{g^\a_{i^n, 0}}\glmpt\grmpt\gcl1\gcl1\gnl
  \gvac2\gev\gnl
  \gend\, \right)=\lambda^2 d_i \nu_{n,n-1}(X_i)=\lambda^2 d_i
  \overline{\nu_n(X_i)}
\end{aligned}
$$
where $N^0_{i^n}=\dim \CC(X_0, X_i^{\o n})$. The third and the
fourth equalities follow from definitions and results in
\cite{NS05}.
\end{proof}
\section{Indicators and the twist of the
center}\label{sec:FS-ribbon} We continue to consider a strict
spherical fusion category $\CC$ over $\BC$. The center $Z(\CC)$ of
$\CC$ is a ribbon category with the twist $\theta$ associated with
the pivotal structure of $\CC$. In this section, we obtain a
formula, in Theorem \ref{t:nu_tr_formula}, for the $n$-th Frobenius
Schur indicator $\nu_n(V)$ of an object $V$ in a spherical fusion
category $\CC$ over $\BC$ in terms of the twist $\theta$ of the
center. The result is a categorical generalization of the formula
for higher indicators of Hopf algebras in \cite[6.4,
Corollary]{KSZ}. By \cite{MugerII03}, $Z(\CC)$ is a modular tensor
category and hence $\theta$ is of finite order (cf. \cite{Vafa88},
\cite{BaKi}). The formula implies that the sequence $\{\nu_n(V)\}_n$
is periodic for each object $V$ and $\nu_n(V)=d(V)$ if $n$ is a
multiple of the order of $\theta$. In addition, if $d(V)$ is
positive for simple $V$, we show in Proposition \ref{t:dD_relation}
that $d(V)=\nu_n(V)$ if and only if $n$ is a multiple of the order
of $\theta$, and for any $n$
 the inequality $|\nu_n(V)| \le d(V)$ holds for all
$V \in \CC$. By a result of \cite[Section 8]{ENO}, there is at most
one pivotal structure on a fusion category such that the pivotal
dimension $d_\ell(V)$ of a simple object $V$ is positive dimensions.
In this case, the pivotal structure is spherical and $d(V)$ is the
Frobenius-Perron dimension of $V$. \sprung

Let $\hat{\Gamma}$ be the set of isomorphism classes of simple
objects of $Z(\CC)$.  For any $(X, e_X) \in \hat{\Gamma}$,
$\theta_{(X,e_X)}=\omega_{(X,e_X)}\id_{(X,e_X)}$ for some root of
unity $\omega_{(X,e_X)} \in \BC$. By M\"uger's results
\cite[Proposition 5.4 and 5.5 ]{MugerII03}, for any $(X,e_X) \in
\hat{\Gamma}$, the element $z_{(X,e_X)} \in \Theta_\CC$ defined by
$$
z_{(X,e_X)}[i,j]=\frac{d(X)}{\lambda d_i} \sum_{\a=1}^{N_i^X}
\left(\,
  \gbeg{4}{7}
  \got1{X_i}\gvac1\got1{X_j} \gnl
  \gcl1\gvac1\gcl1\gnl
  \gdnot{\iota^\a_{i, X}}\glmpt\grmpb\gcl1\gnl
  \gvac1\gbr\gnl
  \gvac1\gcl1\gdnot{\pi^\a_{i, X}}\glmpt\grmpb\gnl
  \gvac1\gcl1\gvac1\gcl1\gnl
  \gvac1\gob1{X_j}\gvac1\gob1{X_i}
  \gend\,\,
\right)
$$
is a primitive central idempotent of $\Theta_\CC$ and
$\sum_{(X,e_X)\in \hat{\Gamma}} z_{(X, e_X)} = \unit_{\Theta_\CC}$,
where
 $\{\iota^\a_{i, X}\}_\a$ is a basis for
$\CC(X_i, X)$  and $\{\pi^\a_{i, X}\}_\a$ its dual basis for $\CC(X,
X_i)$, and $N_i^X=\dim \CC(X_i, X)$. By \cite[Lemma
5.17]{MugerII03}, the element $t$ defined in \eqref{eq:t} can be
written as
\begin{equation}\label{eq:muger}
  t=\sum_{(X,e_X) \in \hat{\Gamma}} \w\inv_{(X,e_X)}
  z_{(X,e_X)}.
\end{equation}
By another result of M\"uger, \cite[Proposition 8.1]{MugerII03}, the
forgetful functor $H: Z(\CC) \rightarrow \CC$ has a two-sided
adjoint $K: \CC \rightarrow Z(\CC)$ such that
\begin{equation}\label{eq:K}
K(Y) \cong \bigoplus_{(X, e_X) \in \hat{\Gamma}} \dim (\CC(X, Y))(X,
e_X)
\end{equation}
for any $Y \in \CC$. Now, we can prove our formula for the
Frobenius-Schur indicators of an object in $\CC$.\sprung

\begin{thm}\label{t:nu_tr_formula}
Let $\CC$ be any spherical fusion category over $\BC$ with pivotal
structure $j$ and let $u$ be the Drinfeld isomorphism of $Z(\CC)$.
For any $V \in \CC$,
$$
\nu_n(V) = \frac{1}{\dim \CC} \pTr(\theta_{K(V)}^n),
$$
where $\theta=u\inv j$ is the twist of $Z(\CC)$ associated with $j$.
\end{thm}
\begin{proof} Since Frobenius-Schur indicators as well as pivotal traces are
invariant under tensor equivalences that preserve the pivotal
structures, we may assume that $\CC$ is strict spherical. In this
case $V^{**}=V$ and $j_V=\id_V$ for all $V \in \CC$. Moreover, since
$\CC$ is semisimple and $K$ preserves direct sums, it suffices to
prove the case when $V=X_i$ for some $i \in \Gamma$. By
\eqref{eq:muger} and Lemma \ref{l:phi_t^n}, we have
$$
\lambda^2 d_i \overline{\nu_n(X_i)} = \phi_i(t^n) = \sum_{(X,e_X)
\in \hat{\Gamma}} \w^{-n}_{(X,e_X)} \phi_i(z_{(X,e_X)}) =
\sum_{(X,e_X) \in \hat{\Gamma}} \w^{-n}_{(X,e_X)}d(X) N_i^X d_i\,,
$$
where $N_i^X = \dim \CC(X_i, X)$. By \cite[Corollary 2.10]{ENO},
$d(V)=d(V^*)$ is real for all $V \in \CC$. Therefore,
\begin{equation}\label{eq:FS-formula1}
\begin{aligned}
\nu_n(X_i)  & =  \overline{\frac{1}{\lambda^2} \sum_{(X,e_X) \in \hat\Gamma} \w^{-n}_{(X,e_X)}d(X) N_i^X} \\
& = \frac{1}{\lambda^2} \sum_{(X,e_X) \in \hat\Gamma} \w^n_{(X,e_X)}d(X) N_i^X \\
&= \frac{1}{\dim \CC} \ptr(\theta^{n}_{K(X_i)})\,. \qedhere
\end{aligned}
\end{equation}
\end{proof}
\begin{rem}
Since
$$
 \ptr(\theta_{K(V^*)}) = \sum_{(X^*,e_{X^*}) \in \hat{\Gamma}}
N_{V^*}^{X^*} \w_{(X^*, e_{X^*})}d(X^*)= \sum_{(X^*,e_{X^*}) \in
\hat{\Gamma}} N_V^X \w_{(X, e_X)}d(X)= \ptr(\theta_{K(V)})
$$
for $V \in \CC$, Theorem \ref{t:nu_tr_formula} implies that
$\nu_n(V)=\nu_n(V^*)$ for all positive integers $n$, which has been
proved in \cite{NS05} using graphical calculus.\sprung
\end{rem}

\begin{rem} \label{rem:nu_formula_v2}
  By Remark \ref{rem:rightcenter}, $\theta'_{(V, e_V)} = \theta\inv_{(V,e_V)}$ for
  object $(V,e_V) \in Z(\CC)$. Then $n$-th Frobenius-Schur indicator
  of $V$ in $\CC$ can be
  rewritten as
  $$
  \nu_n(V) =
  \frac{1}{\dim \CC} \ptr({\theta'}^{-n}_{K(V)})\,.
  $$
  It follows from the equivalence of $Z'(\CC)$ and $\ol Z (\CC)$ that we also have
  $$
  \nu_n(V) =\frac{1}{\dim \CC} \ptr(\bar{\theta}^{-n}_{\ol{K}(V)})
  $$
  where $\ol{K}$ is the two-sided adjoint for the forgetful functor from $\ol Z(\CC)$ to $\CC$ and $\bar{\theta}$
  is the twist of $\ol Z(\CC)$ associated with the pivotal structure of $\CC$.\sprung
\end{rem}
  \begin{rem}
  Since $\theta$ has finite order, the
  sequence $\{\nu_n(V)\}_n$ is periodic for any $V \in \CC$,
  with a period that divides the order of $\theta$.\sprung
  \end{rem}

  The following equation \eqref{eq:dimV} can be also obtained easily by \cite[Proposition 5.4]{ENO}
  using the fact that $\sum_{i \in \Gamma}|d_i|^2 = \dim \CC$. In the following proposition,
  we give another proof for the formula. The special case of the equation for
  $V=I$ is  the class equation in \cite[Proposition 5.7]{ENO}.\sprung

  \begin{prop}\label{t:dD_relation}
    Let $\CC$ be a spherical fusion category over $\BC$ and $\theta$ the twist of $Z(\CC)$
    associated with the pivotal structure of $\CC$. Then
    \begin{equation}\label{eq:dimV}
    d(V)=\frac{1}{\dim \CC}\sum_{(X,e_X) \in \hat{\Gamma}} \dim(\CC(X, V)) d(X)
    \end{equation}
    for all $V \in \CC$. In particular, if $\theta_{K(V)}^n=\id_{K(V)}$, then
    $\nu_n(V)=d(V)$. In addition, if $d_i > 0$ for all $i \in \Gamma$, then the converse also
    holds, and we have
    $$
    |\nu_r(V)| \le d(V)
    $$
     for all positive integers $r$ and $V \in \CC$.
  \end{prop}
  \begin{proof}
    Let $m$ be the order of  $\theta$. Then $\w_{(X,e_X)}^m=1$ for all
    $(X,e_X) \in \hat{\Gamma}$. In particular, $t^m = \unit$ by \eqref{eq:muger}.
    By Lemma \ref{l:phi_t^n}, we have
    $$
    \lambda^2 d_i \overline{\nu_m(X_i)} = \phi_i(\unit)=\lambda^2d_i^2.
    $$
    Since $d_i$ is real, by Theorem \ref{t:nu_tr_formula}, we have
    $$
    d_i = \nu_{m}(X_i) = \frac{1}{\dim\CC} \ptr(\id_{K(X_i)})=\frac{1}{\dim \CC}
\sum_{(X,e_X)\in \hat\Gamma} N_i^X d(X)\,.
$$
The first statement follows directly from the additivity of the
dimension function. In particular, if $\theta_{K(V)}^n=\id_{K(V)}$,
then
$$
\nu_n(V)=\frac{1}{\dim\CC}\ptr(\id_{K({V})})=\frac{1}{\dim \CC}
\sum_{(X,e_X)} {\dim (\CC(X,V))} d(X)=d(V)\,.
$$
If, in addition, $d_i > 0$ for all $i \in \Gamma$, then $d(V)
>0$ for all non-zero $V \in \CC$. Thus, for any positive integer $r$ and for
any $V \in \CC$, we have
$$
|\nu_r(V)|= \left|\frac{1}{\dim \CC}\sum_{(X,e_X) \in \hat{\Gamma}}
\dim (\CC(X, V)) d(X) \w^r_{(X, e_X)}\right| \le \frac{1}{\dim
\CC}\sum_{(X,e_X) \in \hat{\Gamma}} \dim (\CC(X, V)) d(X)=d(V)\,.
$$
Moreover, if $\nu_n(V)=d(V)$, then $\w_{(X,e_X)}^{n}=1$ for any
component $(X,e_X)$ of $K(V)$. Hence, $\theta^n_{K(V)}=\id_{K(V)}$.
\end{proof}
\begin{rem}
Note that
\begin{equation}\label{eq:nu1}
\delta_{i,0} = \nu_1(X_i)={\frac{1}{\dim \CC}}\sum_{(X,e_X) \in
\hat{\Gamma}} N_i^X d(X)\w_{(X,e_X)}
\end{equation}
for $i \in \Gamma$. If $d_i > 0$ for all $i \in \Gamma$,
\eqref{eq:dimV} implies that $\w_{(X,e_X)} =1$ whenever $N_0^X \ne
0$. Thus, the class equation in \cite{ENO} is the special case of
the equation \eqref{eq:nu1} when $i=0$.
\end{rem}

\section{Frobenius-Schur exponent}\label{sec:FSexp}
In this section, we define the {\em Frobenius-Schur exponent} of a
pivotal category over any field $k$, and give an example to
demonstrate the difference between the Frobenius-Schur exponent and
the (quasi)exponent, in the sense of Gelaki and Etingof, of a
spherical fusion category. We prove in Theorem \ref{t:FS_exponent}
that the Frobenius-Schur exponent of a spherical fusion category
$\CC$ over $\BC$ is equal to the order of the twist $\theta$ of the
$Z(\CC)$ associated with the pivotal structure of $\CC$. It will be
shown later in section \ref{s:MTC} that the Frobenius-Schur exponent
of a spherical fusion category is invariant under the center
construction. We then apply Theorems \ref{t:nu_tr_formula} and
\ref{t:FS_exponent} to a semisimple quasi-Hopf algebra $H$ over
${\BC}$ to obtain a formula for the $n$-th indicator of an
$H$-module in terms of the value $\hat{\chi}_V((gu)^{-n})$ of the
character $\hat{\chi}_V$ of the induced module $D(H)\o_H V$, where
$g$ is the trace-element of $H$ (cf. \cite[Section 6]{MN05}) and $u$
is the Drinfeld element of $D(H)$. When $H$ is an ordinary Hopf
algebra, then $g=1$ and the formula reduces to the one in \cite[6.4,
Corollary]{KSZ}. \sprung

As in \cite{NS05},  a pivotal category over a field $k$ is a
$k$-linear pivotal monoidal category with a simple neutral object,
finite-dimensional morphism spaces, and $\End(V)=k$ for all simple
object $V$.\sprung

\begin{defn}
  Let $\CC$ be a pivotal category over any field $k$. The {\em Frobenius-Schur exponent of an object} $V$
  in $\CC$, denoted by $\FSexp(V)$, is defined to be the least positive integer $n$ such that
  $$
  \nu_n(V)=d_{\ell}(V),\,
  $$
  where $d_{\ell}(V)$ and $\nu_n(V)$, respectively, denote the left pivotal dimension and the $n$-th
  Frobenius-Schur indicator of $V$. If such an integer does not exist, we define $\FSexp(V)=\infty$.
  We call $\sup\limits_{V \in \CC}\FSexp(V)$ the {\em Frobenius-Schur exponent of $\CC$} and
  denote it by $\FSexp(\CC)$.\sprung
\end{defn}

 If $H$ is a finite-dimensional semisimple Hopf algebra over
$\BC$, $\FSexp(V)$ is identical to the exponent of $V$ for $V \in
\C{H}$ defined in \cite{KSZ}. Moreover, the results of \cite{KSZ}
show that $\FSexp(H)$ is the same as the exponent of $\C{H}$ in the
sense of Etingof and Gelaki (cf. \cite{EG02} and \cite{Etingof02}).\sprung

\begin{rem} \label{r:fsexp}Let $\CC$ be a pivotal category over a field $k$
with the pivotal structure $j$.
  For a simple object $V\in\CC$,  if $j'$ is another pivotal
  structure, then $j_V$ and $j'_V$ differ by a scalar factor, say
  $j'_V=\lambda j_V$ with $\lambda\in k$. From the definition of
  the left pivotal trace, it is clear that the left pivotal
  dimension $d'_\ell(V)$ computed with respect to $j'$ is
  $d'_\ell(V)=\lambda\inv d_\ell(V)$. Equally, one can read off from the
  definition of the Frobenius-Schur indicators that the indicator
  $\nu'_n(V)$ computed with respect to $j'$ is
  $\nu'_n(V)=\lambda\inv \nu_n(V)$. In particular the
  Frobenius-Schur exponent of a simple object does not depend on
  the choice of a pivotal structure, and neither does the Frobenius-Schur
  exponent of a pivotal fusion category over $k$.\sprung
\end{rem}

Recall from \cite{ENO} that a fusion category $\CC$ over $\BC$ is
pseudo-unitary if it admits a spherical pivotal structure such that
the pivotal dimension $d(V)$ for $V \in \CC$ is positive.\sprung
\begin{prop}
  Let $\CC$ be a pseudo-unitary fusion category over $\BC$. If the
  object $V\in\CC$ contains every simple object of $\CC$, then
  $\FSexp(V)=\FSexp(\CC)$. In particular, if $H$ is a semisimple
  complex quasi-Hopf algebra, then $\FSexp(H)=\FSexp(\C H)$.
\end{prop}
\begin{proof}
  Let $V\cong\sum_{i\in\Gamma}n_iX_i$ and $n=\FSexp(V)$. Then by
  additivity of $\nu_n$
  $$
  \sum_{i \in \Gamma} n_id_i = d(V)=\nu_n(V) = \sum_{i \in \Gamma} n_i \nu_n(X_i)\,.
  $$
  Since $|\nu_n(X_i)| \le d_i$ for $i \in \Gamma$ and $n_i > 0$ by assumption,
  the equality implies that $d_i = \nu_n(X_i)$ for all
  $i \in \Gamma$.
  Now it follows from the additivity of $\nu_n$ that the Frobenius-Schur exponent of $\CC$ is $n$.\sprung

  If $H$ is a semisimple complex quasi-Hopf algebra,
  $\CC=\C{H}$ is a pseudo-unitary fusion category over $\BC$.
  We have $\sum_{i \in \Gamma} d_i X_i \cong H\in\C H$. Therefore, $\FSexp(H)=\FSexp(\CC)$.
\end{proof}

It is reasonable to conjecture that the Frobenius-Schur exponent of
a semisimple quasi-Hopf algebra is identical to its exponent.
However, the following example demonstrates the difference between
the two exponents of a quasi-Hopf algebra and hence $\FSexp(\CC)$ is
different from the exponent of $\CC$ in general.\sprung

\begin{exmp}\label{ex:1}
Let $G=\{1,x\}$ be an abelian group of order 2 and $\w$ a 3-cocycle
of $G$ given by
$$
\w(a,b,c)=\left\{\begin{array}{rl}
-1 &\mbox{if }a=b=c=x;\\
1 & \mbox{otherwise}.
\end{array}\right.
$$
The dual group algebra $\BC[G]^*$ is a well-known Hopf algebra with
the usual multiplication, comultiplication $\Delta$, counit $\e$,
and antipode $S$. Let $\{e(1), e(x)\}$ be the dual basis of $\{1,
x\}$ for $\BC[G]^*$. Define
$$
\phi=\sum_{a,b,c \in G} \w(a,b,c) e(a) \o e(b) \o e(c),\quad
\a=1_{\BC[G]^*},\quad \b=\sum_{a \in G} \w(a,a\inv,a)e(a) =
e(1)-e(x)\,.
$$
Then $H=(\BC[G]^*, \Delta, \e, \phi, \a, \b, S)$ is a quasi-Hopf
algebra and $D(H)=D^\w(G)$. The universal $R$-matrix is given by
$$
R=\sum_{a ,h \in G} e(a) \o 1 \o e(h) \o a\,.
$$
Then
$$
R_{21}R=\sum_{a,b,h,k \in G} (e(a) \o 1)\cdot (e(k) \o b) \o (e(b)
\o 1)\cdot(e(h) \o a)=\sum_{a,b\in G} e(a) \o b \o e(b) \o a\,.
$$
Since
$$
(R_{21}R)^2=\sum_{a,b\in G} \theta_a(b,b)\theta_b(a,a)  e(a) \o b^2
\o e(b) \o a^2 = \sum_{a,b \in G} e(a) \o 1 \o e(b)\o 1=1_{D(H)}\,,
$$
the exponent of $\C{H}$, in the sense of Etingof and Gelaki, is
2.\sprung

  The Frobenius-Schur indicators of the nontrivial
  representation of $H$, on the other hand,
   were already computed in \cite{NS052}. In
  fact $H\cong\BC[G]_x$ is a special case of the constructions
  studied in section 5 of \cite{NS052}, and hence the indicators
  of the nontrivial simple $H$-module $\ol V$ are
  $$
  \nu_n(\ol{V})= \cos\left(\frac{(n-1)\pi}{2}\right)\,.
  $$
  In particular, the Frobenius-Schur exponent of $H$ is 4.\sprung
  \end{exmp}

\begin{thm} \label{t:FS_exponent}
Let $\CC$ be a spherical fusion category over $\BC$. The
Frobenius-Schur exponent of $\CC$ is equal to the order of the twist
$\theta$ of $Z(\CC)$ associated with the pivotal structure of $\CC$.
In particular, $\FSexp(\CC)$ is finite.
\end{thm}
\begin{proof}
Since the Frobenius-Schur exponent as well as the order of $\theta$
are invariant under monoidal equivalences of $\CC$ preserving the
pivotal structure, it suffices to prove the claim in the case when
$\CC$ is a strict spherical category. It follows from
 \eqref{eq:FS-formula1} that
\begin{equation}\label{eq:nu_dim}
 \sum_{i\in
\Gamma} \nu_n(X_i)d_i = \frac{1}{\dim(\CC)}\sum_{(X, e_X) \in
\hat{\Gamma}} \w^n_{(X,e_X)}\sum_{i\in \Gamma}  N_i^Xd_i d(X)
=\frac{1}{\dim(\CC)}\sum_{(X, e_X) \in \hat{\Gamma}} \w^n_{(X,e_X)}
d(X)^2\,.
\end{equation}
If $\nu_n(V)=d(V)$ for all $V \in \CC$, then \eqref{eq:nu_dim}
becomes
$$
\dim(\CC)^2= \sum_{(X, e_X) \in \hat{\Gamma}} \w^n_{(X,e_X)}
d(X)^2\,.
$$
Since $\dim(Z(\CC)) = \dim(\CC)^2$ (cf. \cite{MugerII03}), we have
\begin{equation}\label{eq:dimZ(C)}
   \sum_{(X, e_X) \in \hat{\Gamma}}d(X)^2= \dim(Z(\CC)) = \sum_{(X,
e_X) \in \hat{\Gamma}} \w^n_{(X,e_X)} d(X)^2\,.
\end{equation}
 By
\cite{ENO}, $d(X)$ is real for $(X,e_X) \in \hat{\Gamma}$. The
equation \eqref{eq:dimZ(C)} implies that $\w^n_{(X,e_X)}=1$ for all
$(X,e_X) \in \hat{\Gamma}$. Therefore,
$\theta^n_{(X,e_X)}=\id_{(X,e_X)}$. \sprung

Conversely, suppose that $\theta^n=\id$. By Theorem
\ref{t:nu_tr_formula}, we obtain
$$
\nu_n(V)=\frac{1}{\dim\CC}\ptr(\id_{K(V)})=\frac{1}{\dim \CC}
\sum_{(X,e_X)\in \hat\Gamma} N_i^X d(X)=d(V)
$$
for any $V \in \CC$.
\end{proof}
Since the Drinfeld isomorphism is clearly invariant under
equivalences of braided monoidal categories, we can immediately
conclude:\sprung
\begin{cor}\nmlabel{Corollary}{onlycenter}
  The Frobenius-Schur exponent of a spherical fusion category
  $\CC$ depends only on the equivalence class of the
  spherical braided monoidal category $Z(\CC)$.\sprung
\end{cor}

As it will turn out, the Frobenius-Schur exponent of a spherical
fusion category over $\BC$ is actually invariant under the center
construction. This invariance follows immediately from a formula of
higher Frobenius-Schur indicators for a modular tensor category
which will be derived in section \ref{s:MTC}.\sprung

If the fusion category in question is given as the representation
category of a semisimple quasi-Hopf algebra $H$, then the
characterizations of the Frobenius-Schur indicators and exponent in
terms of the Drinfeld isomorphism in the center turn, of course,
into descriptions in terms of the Drinfeld element of the double of
$H$. The following corollary spells out the details, and the remark
following it gives some more information about the explicit form of
the elements of $D(H)$ involved. Here $D(H)$ stands for a version of
the Drinfeld double construction for quasi-Hopf algebras matching
the left center construction, so that that $\C{D(H)}$ is equivalent
to the center $Z(\C{H})$ of $\C{H}$ as in \cite{Sch:HMDQHA}. See
Remark \ref{comparison} below.\sprung
\begin{cor}\label{c:5.3}
  Let $H$ be a finite-dimensional semisimple quasi-Hopf algebra over
  $\BC$ and $u$ the Drinfeld element of the quantum double $D(H)$ of
  $H$ and $g$ the trace-element of $H$. Then for any simple $H$-module $V$ of $H$,
  $$
  \nu_n(V) = \frac{1}{\dim(H)} \hat{\chi}_V((gu)^{-n})
  $$
  where $\hat{\chi}_V$ is the induced character of the character
  $\chi_V$ of $V$ to $D(H) \o_H V$. Moreover, the Frobenius-Schur exponent of $H$ is equal to
  the order of $gu$.
\end{cor}
\begin{proof}
  Note that the  canonical pivotal structure on $\C{H}$ is given by the formula
  $$
  j_V(x) = g\inv x
  $$
  for any $V \in \C{H}$ and $x \in V$. Moreover, $d(V)=\dim(V)$ for any $V \in
  \CC$. Therefore, the pivotal trace
  of any $f$ in $\End_H(V)$ is identical to the ordinary trace of the linear map
  $f$ and $\C{H}$ is spherical.
  The Drinfeld isomorphism $u_Y$ of $\C{D(H)}$ is given by the
  action of $u$ on the $D(H)$-module $Y$ and the associated twist $\theta$ is given
  by the action of $(gu)\inv$.  Since we always have the
  natural isomorphism
  $$
  \Hom_{D(H)}(D(H) \o_H V, Y) \cong \Hom_H(V, Y),
  $$
  by the uniqueness of adjoint functors, $D(H) \o_H -$ is naturally
  equivalent to $K$. By Theorem \ref{t:nu_tr_formula}, we have
  $$
  \nu_n(V) = \frac{1}{\dim(\C{H})}\Tr_{D(H) \o_H V}((gu)^{-n})
  =\frac{1}{\dim(H)}\hat{\chi}_V((g u)^{-n})\,.
  $$
  The second statement follows immediately from Theorem \ref{t:FS_exponent}.
\end{proof}
\begin{rem}\label{comparison}
  The Drinfeld double $\ol D(H)$ of a finite dimensional quasi-Hopf algebra $H$ is
  usually defined so that $\C{\ol D(H)}=\ol Z (\C{H})$
  (cf. \cite{Kassel} and \cite{Mont93bk} for the Drinfeld double of a Hopf algebra).
  If $H$ is a complex semisimple quasi-Hopf algebra, it follows from Remark \ref{rem:nu_formula_v2}
  that
  $$
  \nu_n(V) = \frac{1}{\dim H} \hat{\chi}_V((gu)^n)
  $$
  for $V \in \C{H}$, where $u$ is the Drinfeld element of $\ol D(H)$ and $g$ is the trace-element
  of $H$. This formula for higher indicators recovers the one in \cite[6.4, Corollary]{KSZ}
  when $H$ is a Hopf algebra.\sprung
\end{rem}
\begin{rem}
  Let $H$ be a quasi-Hopf algebra. By \cite{AC92}, the Drinfeld element $u$
  of $D(H)$ is given by
  $$
  u=S({\phi}^{(-2)} \b S({\phi}^{(-3)}))S(R^{(2)}) \a
  R^{(1)}{\phi}^{(-1)},
  $$
  where ${\phi}^{(-1)} \o {\phi}^{(-2)}\o {\phi}^{(-3)}$ is the
  inverse of the associator of $D(H)$.\sprung

  If $H$ is a semisimple complex quasi-Hopf algebra, one can find
  a formula for $gu$ that does not contain the element $g$
  corresponding to the pivotal structure explicitly. Using the
  expressions for $g$ in \cite[Corollary 8.5]{MN05}, or using
  \cite[Lemma 3.1]{Sch04}, we find that for any $t\in H\o H$
  satisfying $S(t\so 1)\alpha t\so 2=1$ we have
  \begin{align*}
    gu&=gS(\phi\som 2\beta S(\phi\som 3))S(R\so 2)\alpha R\so
    1\phi\som 1\\
    &=gS^2(\phi\som 3)S(\beta)S(\phi\som 2)S(R\so 2)\alpha R\so
    1\phi\som 1\\
    &=\phi\som 3gS(\beta)S(\phi\som 2)S(R\so 2)\alpha R\so 1\phi\som
    1\\
    &=\phi\som 3\Lambda\sw 2t\so 2S(\Lambda\sw 1t\so 1)S(\phi\som
    2)S(R\so 1)\alpha R\so 1\phi\som 1.
  \end{align*}
  Possible choices are $t=p_L$ or $t=p_R$.\sprung
\end{rem}

We have already computed the Frobenius-Schur exponent of the
nontrivial two-dimensional quasi-Hopf algebra above, and seen that
it differs from the exponent in the sense of Etingof. We shall redo
this example now using the new characterization of the
Frobenius-Schur exponent as the order of $gu$:\sprung

\begin{exmp}
Let $H$ be the quasi-Hopf algebra given in Example \ref{ex:1}.
 Since $D^\w(G)$ is commutative
and the antipode of $D^\w(G)$ is the identity map,  the general
formula for $u$ from \cite{AC92} specializes to
$$
u=R^{(2)}R^{(1)} = \sum_{a \in G} e(a) \o a = e(1)\o 1 + e(x) \o x.
$$
Direct computation shows that $\ord(u)=4$. By \cite{MN05}, the
trace-element $g$ of $D(H)$ is given by
$$
g=\sum_{a \in G} \w(a,a\inv,a)e(a) \o 1 = (e(1)-e(x))\o 1\,.
$$
Therefore, $\ord(g)=2$. By the commutativity of $D^\w(G)$ again,
$$
\ord(gu) =4.
$$
Therefore, as we have already seen in Example \ref{ex:1}, the
Frobenius-Schur exponent of $H$ is 4.\sprung
\end{exmp}

Since the canonical pivotal structure of the module category over a
semisimple quasi-Hopf algebra is preserved by every monoidal
equivalence, we can deduce some invariance properties of the
Frobenius-Schur exponent:\sprung
\begin{prop}\nmlabel{Proposition}{onlydouble}
Let $H,H'$ be complex semisimple quasi-Hopf algebra.
\begin{enumerate}
  \item\label{DHinv} The Frobenius-Schur exponent of $H$ depends only on the
  gauge equivalence class of the double $D(H)$ as a quasitriangular quasi-Hopf
  algebra, i.e. $\FSexp(H)=\FSexp(H')$ provided $D(H)$ and $D(H')$
  are gauge equivalent quasitriangular quasi-Hopf
  algebras.
  \item\label{tpexp} $\FSexp(H\o H')=\lcm(\FSexp(H),\FSexp(H'))$.
  \item\label{opexp} $\FSexp(H^\op)=\FSexp(H^\cop)=\FSexp(H)$.
  \item\label{DHexp} $\FSexp(D(H))=\FSexp(H)$.
\end{enumerate}
\end{prop}
\begin{proof}
  Since the pivotal structure is preserved under any monoidal
  equivalence between module categories of semisimple complex quasi-Hopf
  algebras, \eqref{DHinv} is a direct consequence of
  \nmref{onlycenter}. Statement \eqref{tpexp} can be verified
  directly from the definition of the Frobenius-Schur exponent,
  since one easily sees $\nu_n(V\o W)=\nu_n(V)\nu_n(W)$ for
  $V\in\C H$ and $W\in\C{H'}$. As for \eqref{opexp}, it suffices
  to treat $H^\cop$, since $H^\op$ and $H^\cop$ are gauge
  equivalent through the antipode. Now $\C{H^\cop}=(\C H)^\sym$ is
  the category $\C H$ with the opposite tensor product. For
  $n=\FSexp(H)$ we have
  $\dim(V)=\ol{\dim(V)}=\ol{\nu_n(V)}=\nu_{n,n-1}(V)=\nu_n(V^\sym)$ by
  \cite[Theorem 5.1, Lemma 5.2]{NS05}, where $V^\sym$ denotes the module
  $V$ considered as an object of $(\C H)^\sym$. Thus
  $\FSexp(H)$ divides $\FSexp(H^\cop)$, and by symmetry we are
  done. Finally \eqref{DHexp} follows since $D(D(H))$ is gauge equivalent
  to $D(H\o H^\op)$. This is can be rephrased and proved entirely in categorical
  terms (see the remark below). For quasi-Hopf algebras we can
  argue as follows, without using semisimplicity: By
  \cite{Sch:HMDQHA} the category of modules over the double $D(H)$
  is isomorphic to the monoidal category $^H\CC^H$ of $H$-$H$-bicomodules in
  the monoidal category $\CC=\C H$-$H$ of $H$-bimodules. By \cite{Schauenburg01}
  the center of this bicomodule category is equivalent to the
  center of the underlying category, and so
  $\C{D(H\o H^\op)}\cong Z(\C{H\o H^\op})\cong Z(\CC)\cong Z({^H\CC^H})\cong Z(\C{D(H)})\cong\C{D(D(H))}$
  as braided monoidal categories.
\end{proof}
\begin{rem}
For any fusion category $\CC$, we have a braided monoidal category
equivalence $Z(Z(\CC))\cong Z(\CC\boxtimes\CC^\sym)$ by
\cite[section 7]{MugerII03}. Thus at least in the semisimple case we
need here the result on the double of the double of a quasi-Hopf
algebra above has a categorical version using a suitable tensor
product of categories. However, we did not verify that the
equivalence preserves the pivotal structures in this case, so we
cannot draw the desired conclusion on the Frobenius-Schur exponent
of the center for general spherical fusion categories. We will
arrive at that result with an entirely different proof in section
\ref{s:MTC}.\sprung
\end{rem}

We close this section with\sprung
\begin{cor}
  Let $\CC$ be a spherical fusion category over $\BC$ with pivotal structure $j$.
  Then for any simple $V \in \CC$, $I$ is a summand of $V^{\o n}$, where $n=\FSexp(V), \FSexp(\CC)$.
  Moreover, for any pivotal structure $j'$ on $\CC$, $j\inv j'$ is a finite
  order monoidal automorphism of the identity functor of $\CC$, and
  $\ord(j\inv j') \mid \FSexp(\CC)$.
\end{cor}
\begin{proof}
  By \cite{NS05},  $\nu_r(V)$ is the ordinary trace of an $\BC$-linear automorphism
  $E_V^{(r)}$ on $\CC(I, V^{\o r})$.  Thus, if $\nu_r(V) \ne 0$,
  then $\CC(I, V^{\o r}) \ne 0$. By \cite[Theorem 2.3]{ENO}, $d(V)\ne 0$ for any simple object $V$,
  and $\nu_n(V)= d(V)$.
  Therefore, $\CC(I, V^{\o n}) \ne 0$. Let $\lambda = j\inv j'$. Then the component $\lambda_V$ is a
  scalar, and $\lambda_I=1$ (cf. Remark \ref{r:piv_dim}).
  Since $I$ is a summand of $V^{\o n}$, we have $\lambda_V^n =1$ and
  hence $\ord(\lambda_V) \mid \FSexp(\CC)$. Since $\ord(\lambda)=\lcm \ord(\lambda_V)$,
  where $V$ runs through a complete set of non-isomorphic simple
  objects of $\CC$, the divisibility $\ord(\lambda) \mid
  \FSexp(\CC)$ follows.
\end{proof}

The fact that $I$ is a direct summand of some tensor power of every
simple $V$ generalizes a result from \cite[Section 4]{KSZ}; in their
terminology (introduced for modules over Hopf algebras), any simple
$V$ has finite order.

\section{Etingof's exponent vs. Frobenius-Schur
exponent}\label{sec:versus} We have discussed already that the
Frobenius-Schur exponent of a quasi-Hopf algebra can differ from its
exponent as defined by Etingof. In the example, the difference
amounts to a factor $2$. The main result of this section implies
that this is the most general discrepancy between the two notions
that can occur. In particular, results about the Frobenius-Schur
exponent have implications for the exponent in the sense of
Etingof.\sprung

Let $\CC$ be a  ribbon fusion category over $\BC$ with the twist
$\theta$ and braiding $c$. Let $M_{V,W}=c_{W,V}\circ c_{V,W}$ for
any $V, W \in \CC$. By \cite{Etingof02}, $\theta$ has finite order,
and $M$ has finite order with $\ord(M)\mid\ord(\theta)$.\sprung
\begin{prop}\nmlabel{Proposition}{atmosttwice}
  $\ord(\theta)=\ord(M)$ or $2\ord(M)$.
\end{prop}
\begin{proof}
Let $X_1, \cdots, X_l$ be the complete set of isomorphism classes of
simple objects of $\CC$, and let $\theta_{X_i}=\w_i \id_{X_i}$ for
some root of unity $\w_i$. Then
$$
X_i^*= X_{\bar{i}}
$$
for some $\bar{i}$. The ribbon structure $\theta$ implies that
$\w_i=\w_{\bar{i}}$ and $\CC$ is spherical. Hence,
$d_i=d_{\bar{i}}$.  Let $n=\ord(M)$. Then the equality
$$
\theta_{V\o W} = (\theta_V \o \theta_W)M_{V, W}
$$
implies that
$$
\theta^n_{V\o W} = (\theta^n_V \o \theta^n_W)
$$
for all $V, W \in \CC$. In particular, we have
\begin{equation}\label{eq:ribbon_eq2}
\theta^n_{X_i\o X_j} = (\theta^n_{X_i} \o \theta^n_{X_j})
\end{equation}
for all $i, j \in \{1,\dots, l\}$. Let $N_{ij}^k=\dim (\CC(X_k, X_i
\o X_j))$. Taking trace on both sides of \eqref{eq:ribbon_eq2}, we
have
$$
\sum_{k \in \Gamma}\w_k^n d_k N_{ij}^k = \w_i^n\w_j^n d_id_j
$$
where $d_i$ is the pivotal dimension of $X_i$. Let
$$
v_i=\w_i^n d_i, \quad \mathbf{v} =\left[\begin{array}{c}
  v_1\\
  \vdots \\
  v_l\\
\end{array}\right] ,\quad N_i = [N^b_{ia}]_{a,b}\,.
$$
Then we have
$$
N_i\mathbf{v}=v_i \mathbf{v}, \quad N_{\bar{i}}=N^t\,.
$$
Taking the complex transposition of the first equation, we have
$$
\bar{\mathbf v}^t N_{\bar{i}} = \bar{v}_i \bar{\mathbf v}^t
$$
and hence
$$
\bar{v}_i  \bar{\mathbf v}^t\mathbf{v}=  \bar{\mathbf v}^t
N_{\bar{i}} \mathbf{v} =
 v_{\bar{i}}\bar{\mathbf v}^t \mathbf{v}\,.
$$
Since $\CC$ is spherical and $\mathbf{v}$ is a non-zero complex
vector, we obtain
$$
\bar{\w}_i^n d_i = \bar{v}_i= v_{\bar{i}} = \w^n_i d_{\bar{i}}=
\w^n_i d_i\,.
$$
Since $d_i \ne 0$, $\bar{\w}^n_i=\w^n_i$. Therefore, $\w^n_i = \pm
1$ and so $\w^{2n}_i=1$ for all $i$. Equivalently,
$\theta^{2n}=\id$.
\end{proof}
\begin{cor}\label{c:6.2}
  For any semisimple quasi-Hopf algebra $H$ over $\BC$ we have
  $\FSexp(H)=\exp(H)$ or $\FSexp(H)=2\exp(H)$. \qed \sprung
\end{cor}

 We close this section with the following lemma which will be used in section \ref{sec:Cauchy}
 to prove the Cauchy Theorem
 for semisimple quasi-Hopf algebras.\sprung
\begin{lem}\nmlabel{Lemma}{l:change_twist}
Let $\CC$ be a ribbon fusion category over $\BC$ with the braiding
$c$, the pivotal structure $j$, and the twist $\theta$. Let $M_{V,W}
= c_{W,V}c_{V,W}$ for $V,W \in \CC$. If $\ord(\theta) = 2 \ord (M)$
and $\ord(M)$ is odd, then there exists a spherical pivotal
structure $\hat{j}$ on $\CC$ such that the order of the twist
$\hat\theta$ associated with $\hat{j}$ is equal to the order of $M$,
and $\hat{d}(V)=\pm d(V)$ where $\hat{d}(V)$ and $d(V)$ denote the
dimensions of $V$ computed with the pivotal structures $\hat{j}$ and
$j$, respectively. In addition, if $\CC$ is a MTC, then $\CC$ with
the twist $\hat{\theta}$ is also a MTC.
\end{lem}
\begin{proof}
  Without loss of generality, we may assume that $\CC$ is strict
  pivotal. Let $N=\ord(M)$ and $\hat{\theta}=\theta^{N+1}$. Since
  $N$ is odd, $\gcd(N+1, 2N)=2$ and so $\ord(\hat\theta) = N$.
  Moreover, $\theta^N_{V\o W} = \theta^N_V \o \theta^N_W$ for any
  $V, W \in \CC$. Thus, $\hat\theta$ is a twist and $\CC$ is a ribbon category
  with respect to $\hat\theta$. Let $\hat{j}=u\hat\theta$ be the
  spherical pivotal structure on $\CC$ associated with $\hat\theta$,
  where $u$ is the Drinfeld isomorphism of the braiding $c$. Since
  $\CC$ is strict pivotal, $u=\theta\inv$ and so $\hat{j}=\theta^N$.
  Let $X_i$, $i \in \Gamma$, be a complete set of non-isomorphic
simple objects of $\CC$. Then, for $i \in \Gamma$,
$\theta^N_{X_i}=\pm \id_{X_i}$ and
  $$
  \hat{d}_i=\hat{d}(X_i)=
  \def\objectstyle{\scriptstyle}
\xy
  (0,3)*{}="o", "o"+(8,0)="l1",
  {"o";"l1" **\crv{"o"+(0,6)& "l1"+(0,6)}},
  "o"+(0,-2.7)*+[F]{\theta_{X_i}^N},
  "o"+(0,-5.5)="o1",
  "o1"+(8,0)="l2",
  {"o1";"l2" **\crv{"o1"+(0,-6)& "l2"+(0,-6)}},
  {"l1"\ar@{-}"l2"}
\endxy
  = \pm d(X_i) = \pm d_i\,.
  $$
  In particular, $\sum_{i \in \Gamma} {\hat{d}}_i^2 = \sum_{i \in
  \Gamma} d_i^2=\dim \CC$.\sprung

  Suppose, in addition, $\CC$ is a MTC with respect to the twist
  $\theta$. Let $\Gamma^\pm=\{i \in \Gamma\mid \theta^N_{X_i} = \pm \id_{X_i}\}$.
   Let $[s_{ij}]_{i,j\in \Gamma}$ be the
   $S$-matrix of this MTC. Now, we compute the $S$-matrix
$[\hat{s}_{ij}]$ of the ribbon category $\CC$ with respect to the
twist $\hat{\theta}$. Note that
$$
  \theta^N_{X^*_i \o X_j} = \theta^N_{X^*_i} \o \theta^N_{X_j} =
  \left\{\begin{array}{rl}
  \id_{X^*_i \o X_j} & \text{if } i, j \in \Gamma^+ \text{ or } i,j
  \in \Gamma^-,\\
  &\\
  -\id_{X^*_i \o X_j} & \text{otherwise}.
\end{array}
\right.
$$
Thus
$$
\hat{s}_{ij} = \frac{1}{\sqrt{\dim \CC}} \cdot
\def\objectstyle{\scriptstyle}
\xy
  (0,8)*{}="o", "o"+(15,0)="l1",
  {"o";"l1" **\crv{"o"+(0,8)& "l1"+(0,8)}},
  "o"+(0,-2.8)*+[F]{M_{X_i^*, X_j}},
  {"o"+(0,-5.6)\ar@{-}"o"+(0,-10)="o1"},
  "o1"+(0,-3)*+[F]{\theta_{X_i^*\o X_j}^N},
  "o1"+(0,-6.2)="o2", "o2"+(15,0)="l2",
  {"o2";"l2" **\crv{"o2"+(0,-8)& "l2"+(0,-8)}},
  {"l1"\ar@{-}"l2"}
\endxy
= \left\{
\begin{array}{rl}
  s_{ij} & \text{if } i,j \in \Gamma^+ \text{ or } i,j \in
  \Gamma^-,\\
  -s_{ij} & \text{otherwise}.
\end{array}
\right.
$$
Therefore, $[s_{ij}]$ and $[\hat{s}_{ij}]$ are conjugate matrices.
Hence the matrix $[\hat{s}_{ij}]$ is non-singular and $\CC$ is a MTC
with respect to $\hat\theta$.
\end{proof}

\section{Frobenius-Schur indicator formula for Modular tensor
Categories}\label{s:MTC}
 In \cite{Bantay97}, Bantay has defined a
scalar, called the (2nd)Frobenius-Schur indicator, by the formula
\begin{equation}\label{eq:bantay_2_formula}
  \sum_{i,j  \in \Gamma} N^k_{ij} s_{0i} s_{0j}
\left(\frac{\w_i}{\w_j} \right)^2
\end{equation}
for each simple object $X_k$ of a modular tensor category $\CC$
(MTC), where $[s_{ab}]$ denotes the $S$-matrix of $\CC$ and
$N_{ij}^k=\dim (\CC(X_i \o X_j, X_k))$. It is shown in the paper
that the value of the expression can only be 0, 1, or $-1$. \sprung

In this section, we will derive a formula (Theorem
\ref{t:bantay_n_formula}) for the $n$-th indicator $\nu_n(V)$ of a
simple object $V$ in a modular tensor category $\CC$ using the
definition of higher indicators introduced in \cite{NS05}. Bantay's
formula \eqref{eq:bantay_2_formula} is recovered as the special case
$n=2$. This also shows that Bantay's notion of Frobenius-Schur
indicator is the 2nd FS-indicator in our sense.\sprung

As an immediate consequence of the formula for $\nu_n(V)$, we show
that $\FSexp(\CC)$ is equal to the order of the twist $\theta$ of
$\CC$ and that the Frobenius-Schur exponent of a spherical fusion
category is invariant under the center construction.\sprung

We continue to use the notation introduced in section
\ref{sec:tube}. Let $\CC$ be a {\em strict} modular tensor category,
i.e. a MTC whose underlying spherical category is strict. Suppose
that $X_i$, $i \in \Gamma$, is a complete set of non-isomorphic
simple objects of $\CC$. We define $\w_i \in \BC$, $i \in \Gamma$,
by the equation
$$
\theta_{X_i} =\w_i \id_{X_i}.
$$
 Since the twist $\theta$ of $\CC$
is of finite order, $\w_i$ is a root of unity.\sprung

By \cite{NS05}, the $(n,k)$-th Frobenius-Schur indicator
$\nu_{n,k}(V)$ of the object $V \in \CC$ is the ordinary trace of
the linear operator $\left(E_V^{(n)}\right)^k : \CC(I, V^{\o n})
\rightarrow \CC(I, V^{\o n})$ defined by
\[
\left(E_V^{(n)}\right)^k \colon\quad \xy
(0,0)*+{\includegraphics{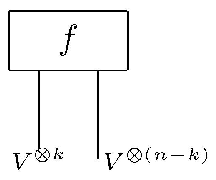}}
\endxy
\quad \mapsto \quad \xy (0,-2)*+{\includegraphics{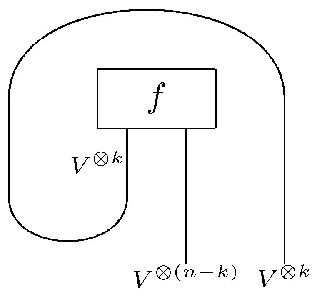}}
\endxy
\quad = \quad \xy (0,-2)*+{\includegraphics{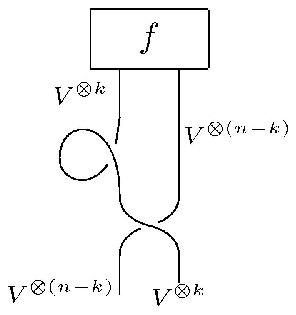}}
\endxy\,.
\]
Let $C_V^{(m,n)} \in \CC(V^{\o(m+n)}, V^{\o(m+n)})$ denote the map
$$
\xy (0,0)*+{\includegraphics{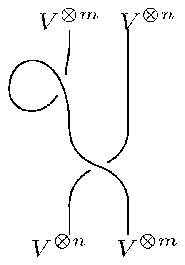}}
\endxy\,.
$$
In particular,
\begin{equation}\label{eq:c(1,n-1)}
 C_{X_j}^{(1, n-1)} = \w_j\cdot c_{X_j, X_j^{\o (n-1)}}\,.
\end{equation}
 For any positive integers $n, k$, let $\{p_\a\}$ be a
basis for $\CC(I, V^{\o n})$, and $\{q_\a\}$ the dual basis for
$\CC(V^{\o n}, I)$. Then we have
$$
\Tr\left(\left(E_V^{(n)}\right)^k\right) = \sum_\a \xy
(0,0)*{\includegraphics{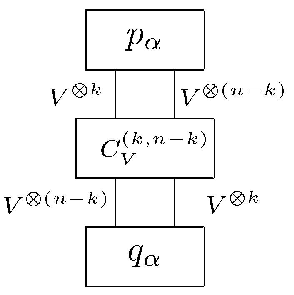}}
\endxy
\quad = \quad \sum_{\a}\quad \xy (0,0)*{\includegraphics{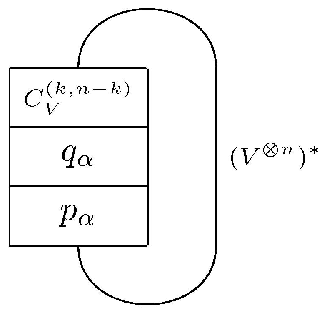}}
\endxy\,.
$$
Let $e_0(V^{\o n})=\sum_\a p_\a q_\a$. Then we have
\begin{equation}\label{eq:MTC_nu}
\nu_{n,k}(V)=\quad \xy (0,0)*{\includegraphics{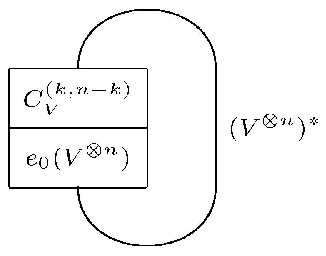}}\endxy\,.
\end{equation}
Note that $e_0(V^{\o n})$ is independent of the choices of the basis
$\{p_\a\}$ for $\CC(I,V^{\o n})$. In general, the map $e_0(W)$, for
any object $W$ of $\CC$, is the idempotent of ${\CC(W, W)}$ given by
\begin{equation}\label{eq:e_0}
e_0(W)=\iota \pi,
\end{equation}
where $\pi: W \rightarrow W\triv$ and $\iota: W\triv \rightarrow W$
are the epimorphism and the monomorphism associated with the summand
$W\triv$ of $W$.\sprung

\begin{lem}\label{l:e0_formula1} For any $W \in \CC$, we have
  $$
  e_0(W) = \frac{1}{\dim(\CC)}\sum_{i \in \Gamma}d_i\,\,
  \xy (0,0)*{\includegraphics{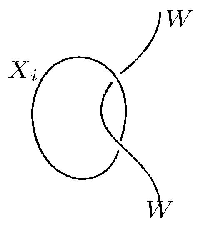}}\endxy\,.
  $$
\end{lem}
\begin{proof}
The statement follows directly from \eqref{eq:e_0} and
\cite[Corollary 3.1.11]{BaKi}.
\end{proof}
By Lemma \ref{l:e0_formula1} and \eqref{eq:MTC_nu}, we have
\begin{equation}\label{eq:MTC_nu2}
 \nu_{n,k}(V) =
  \frac{1}{\dim(\CC)}\sum_{i \in \Gamma}d_i\quad
  \left(
   \xy (0,0)*{\includegraphics{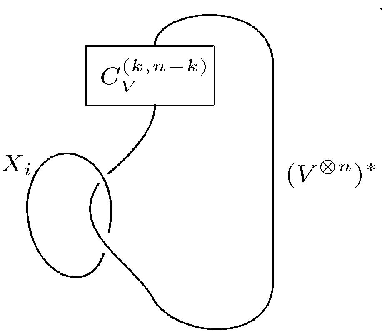}}\endxy
  \right)\,.
\end{equation}
Since $\nu_n(V)=\nu_{n,1}(V)$, the following lemma follows
immediately from \eqref{eq:c(1,n-1)}.

\begin{lem} \label{l:nu1}
For any $j \in \Gamma$ and integer $n \ge {1} $, we have
$$
\nu_{n}(X_j) =
  \frac{1}{\dim(\CC)}\sum_{i \in \Gamma}d_i \w_j\quad
  \left(\xy (0,0)*{\includegraphics{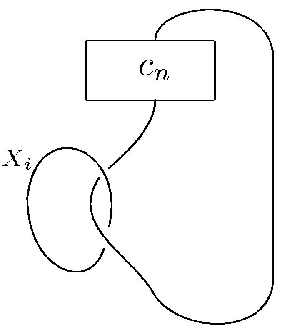}}
\endxy\,\, \right),
$$
where $c_n$ denotes the map $c_{X_j, X_j^{\o (n-1)}}$. \qed\sprung
\end{lem}

Let $V \in \CC$. We define $F: \CC(V\o X_j, V\o X_j)\rightarrow
\CC(V\o X_j, V\o X_j)$ by
\begin{equation}\label{eq:F1}
F\colon \xy (0,0)*{\includegraphics{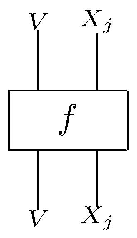}}
\endxy
\quad\mapsto \quad \xy (0,0)*{\includegraphics{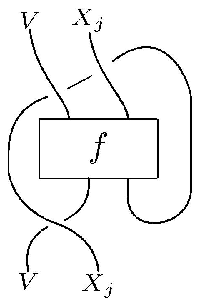}}
\endxy
\quad =  \xy (0,0)*{\includegraphics{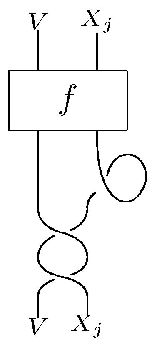}}
\endxy\,.
\end{equation}
Let $m_n=c_{X_j^{\o n}, V} \circ c_{V, X_j^{\o n}}$, and
$c_n=c_{X_j, X_j^{\o (n-1)}}$ for any integer $n \ge 1$. Then,
\eqref{eq:F1} says that
$$
F(f)=\w_j\cdot m_1 \circ f
$$
for $f \in \CC(V\o X_j, V\o X_j)$. In particular,
\begin{equation}\label{eq:F2}
F^{n-1}(m_1) = \w_j^{n-1} \cdot m_1^n\,,
\end{equation}
where $F^0$ denotes the identity map.\sprung
\begin{lem}\label{l:t1}
For any integer $n \ge 0$, we have
  $$
  F^n(m_1)=\quad \xy (0,0)*{\includegraphics{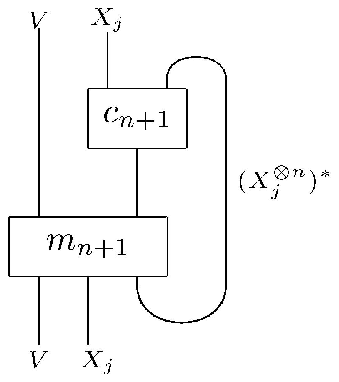}} \endxy\,.
  $$
\end{lem}
\begin{proof}
  The equality is obviously true for $n=0$ as $c_1=\id_{X_j}$. Assume the equation holds for some integer
  $n\ge 0$. Then
 $$
  F^{n+1}(m_1)
 =\quad
 \xy (0,0)*{\includegraphics{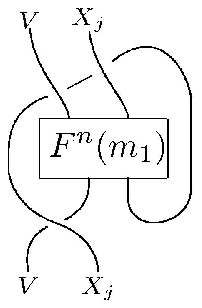}}\endxy
\quad=\quad \xy (0,0)*{\includegraphics{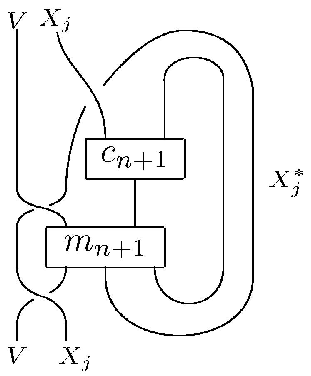}}\endxy\,.
  $$
  Since
  $$
  c_{n+2}=\quad
  \xy (0,0)*{\includegraphics{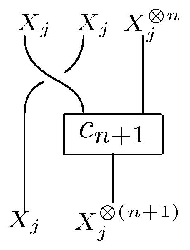}}\endxy
  \quad \mbox{and}\quad m_{n+2} =\quad
  \xy (0,0)*{\includegraphics{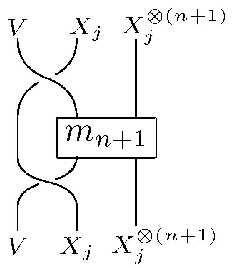}}\endxy\,,
  $$
  we have
    $$
  F^{n+1}(m_1)=\quad
\xy (0,0)*{\includegraphics{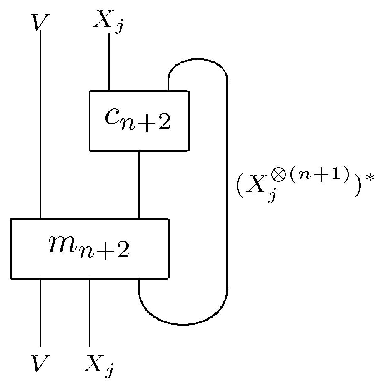}}\endxy\,,
$$
and this completes the proof.
\end{proof}
\begin{lem}\label{l:t2}
  For any integer $n \ge {1}$,
$$
\xy (0,0)*{\includegraphics{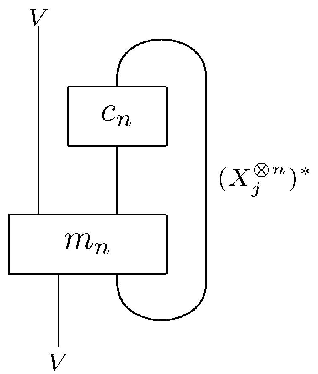}}\endxy
\quad = \quad \w_j^{n-1}\cdot
\xy (0,0)*{\includegraphics{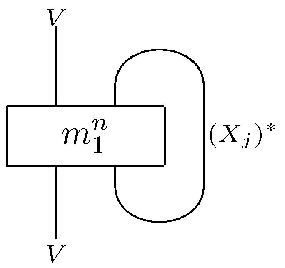}}\endxy
\,.
  $$
\end{lem}
\begin{proof}
  The equality follows immediately from Lemma \ref{l:t1} and \eqref{eq:F2}.
\end{proof}
\begin{thm}\label{t:bantay_n_formula}
  For any $j \in \Gamma$ and positive integer $n$, we have
  $$
  \nu_n(X_j)= \frac{1}{\dim(\CC)}\sum_{i, k \in \Gamma}N_{ik}^j d_i d_k  \left(\frac{\w_i}{\w_k}\right)^n,
  $$
  where $N_{ik}^j=\dim \CC(X_i \o X_k, X_j)$.
  In particular, if $N=\exp(\CC)$, then
  $$
  \nu_N(X_j)= \w_j^N d_j = \pm d_j\,.
  $$
\end{thm}
\begin{proof}
By Lemmas \ref{l:nu1} and \ref{l:t2},
\begin{equation}\label{eq:m1}
\nu_n(X_j)= \frac{1}{\dim(\CC)}\sum_{i \in \Gamma} d_i \w_j^n \cdot
\def\objectstyle{\scriptstyle}
\xy (0,0)="ctext"; "ctext"*{\displaystyle m_1^n}; "ctext"+(-6,
-3)*{}; "ctext"+(6,3)*{} **\frm{-};
"ctext"+(-2,5)="c2";"ctext"+(-2,3)*{} **\dir{-};
"ctext"+(2,3)*{};"ctext"+(2,5)*{}="b1" **\dir{-} ;
"ctext"+(-2,-3)*{};"ctext"+(-2,-5)="d2" **\dir{-};
"ctext"+(2,-3)*{};"ctext"+(2,-5)*{}="a1" **\dir{-} ;
"a1";"a1"+(8,0)="a2" **\crv{"a1"+(0,-5) &
"a2"+(0,-5)}?(0.8)+(3.5,0)*{X_j^*}; "b1";"b1"+(8,0)="b2"
**\crv{"b1"+(0,5) & "b2"+(0,5)}; "b2";"a2" **\dir{-} ;
"c2"-(8,0)="c1"; "c2" **\crv{"c1"+(0,5) & "c2"+(0,5)};
"d2"-(8,0)="d1"; "d2" **\crv{"d1"+(0,-5) & "d2"+(0,-5)}; "c1";"d1"
**\dir{-} ?(0.2)+(-2,0)*{X_i};
\endxy
= \frac{1}{\dim(\CC)}\sum_{i \in \Gamma} d_i \w_j^n \cdot \xy
(0,0)="ctext"; "ctext"*{\displaystyle m_1^n}; "ctext"+(-4, -3)*{};
"ctext"+(4,3)*{} **\frm{-}; "ctext"+(0,3)*{};"ctext"+(0,5)*{}="b1"
**\dir{-} ; "ctext"+(0,-3)*{};"ctext"+(0,-5)*{}="a1" **\dir{-} ;
"a1";"a1"+(8,0)="a2" **\crv{"a1"+(0,-5) & "a2"+(0,-5)};
"b1";"b1"+(8,0)="b2" **\crv{"b1"+(0,5) &
"b2"+(0,5)}?(0,3)+(-6,0)*{X_i^*\o X_j}; "b2";"a2" **\dir{-};
\endxy\,.
\end{equation}
Since $\w_i \w_j m_1 = \theta_{X_i^* \o X_j}$, we have
$$
m_1^n = \frac{1}{(\w_i\w_j)^n}\bigoplus_{k \in
\Gamma}\left(\theta^n_{X_k}\right)^{\oplus N_{\bar{i} j}^k}\,.
$$
Therefore,
$$
\nu_n(X_j)=\frac{1}{\dim(\CC)}\sum_{i,k \in \Gamma} d_i
\frac{\w_k^n}{\w_i^n}N_{\bar{i}j}^k d_k
=\frac{1}{\dim(\CC)}\sum_{i,k \in \Gamma} d_i d_k N_{ik}^j
\frac{\w_i^n}{\w_k^n}\,.
$$
In particular, if $N=\exp(\CC)$, then $m_1^N=\id$. By Proposition
\ref{atmosttwice}, $\w_j^N=\pm 1$. It follows from \eqref{eq:m1}
that
$$
\nu_N(X_j)= \frac{1}{\dim(\CC)}\sum_{i \in \Gamma} d_i d_{\bar{i}}
d_j \w_j^{N} =\w_j^N d_j \,.\qedhere
$$
\end{proof}
\begin{rem} Theorem \ref{t:bantay_n_formula} implies that the Frobenius-Schur indicators of the objects in a MTC are
real. The formula for the $n$-th indicators an be rewritten in term
of the modular data of $\CC$ by the Verlinde formula (cf.
\cite{BaKi}), namely
$$
N_{ik}^j= \sum_{r \in \Gamma} \frac{s_{ir}s_{kr}s_{\bar{j} r}}{s_{0
r}}\,,
$$
where the $S$-matrix of $\CC$ is given by
$$
s_{ij}=\frac{1}{\sqrt{\dim(\CC)}}
\def\objectstyle{\scriptstyle}
\xy (0,0)="ctext"; "ctext"*{\displaystyle m_1}; "ctext"+(-4, -3)*{};
"ctext"+(4,3)*{} **\frm{-}; "ctext"+(0,3)*{};"ctext"+(0,5)*{}="b1"
**\dir{-} ; "ctext"+(0,-3)*{};"ctext"+(0,-5)*{}="a1" **\dir{-} ;
"a1";"a1"+(8,0)="a2" **\crv{"a1"+(0,-5) & "a2"+(0,-5)};
"b1";"b1"+(8,0)="b2" **\crv{"b1"+(0,5) &
"b2"+(0,5)}?(0,3)+(-6,0)*{X_i^*\o X_j}; "b2";"a2" **\dir{-};
\endxy\quad.
$$
In particular, $s_{0i}=d_i/\sqrt{\dim(\CC)}$ for all $i \in \Gamma$.
Hence
$$
\nu_n(X_j)=\sum_{i,k \in \Gamma} N_{ik}^j s_{0i} s_{0k}
\frac{\w_i^n}{\w_k^n}\,.
$$
For $n=2$, this recovers Bantay's formula for the (degree 2)
Frobenius-Schur indicator in conformal field theory (cf.
\cite{Bantay97}).\sprung
\end{rem}

\begin{thm}\label{t:o(theta)}
  Let $\CC$ be a MTC with the twist $\theta$. Then $\FSexp(\CC)=\ord(\theta)$.
\end{thm}
\begin{proof}
  If $\theta^n=\id$,  then
  $$
  \nu_n(X_j) = \frac{1}{\dim(\CC)}\sum_{i,k} N_{ik}^j d_id_k =
  \frac{1}{\dim(\CC)}\sum_i d_i d_{\bar{i}} d_j = d_j
  $$
  for all $j \in \Gamma$. Conversely, if $\nu_n(X_j)=d_j$ for all $j
  \in \Gamma$, then
  $$
  \begin{aligned}
  \dim(\CC)& =\sum_{j} d_j^2 =
  \frac{1}{\dim(\CC)}\sum_{i,k,j \in \Gamma} d_i d_k N_{ik}^j d_j
  \frac{\w_i^n}{\w_k^n} \\
  &=
  \frac{1}{\dim(\CC)}\sum_{i,k \in \Gamma} d_i^2 d_k^2
  \frac{\w_i^n}{\w_k^n}=
   \frac{1}{\dim(\CC)}\left|\sum_{i \in \Gamma} d_i^2
   \w_i^n\right|^2\,.
  \end{aligned}
  $$
  Hence, we have
  $$
  \sum_{i \in \Gamma} d_i^2 =\dim(\CC) = \left|\sum_{i \in \Gamma} d_i^2
   \w_i^n\right|  \,.
  $$
  The equalities imply that $\w_i^n$ are all identical for $i \in
  \Gamma$. Since $\w_0^n=1$, $\w_i^n=1$ for all $i$ and so
  $\theta^n=\id$.
\end{proof}
\begin{cor}\label{c:invariance}
  Let $\CC$ be a spherical fusion category over $\BC$. Then
  $\FSexp(\CC)=\FSexp(Z(\CC))$.
\end{cor}
\begin{proof}
 By  \cite{MugerII03}, $Z(\CC)$ is a modular tensor category and hence
  $\FSexp(Z(\CC))=\ord(\theta)$ by Theorem \ref{t:o(theta)}, where $\theta$ is the twist
  of $Z(\CC)$ associated with the pivotal structure of $\CC$.
  By Theorem \ref{t:FS_exponent}, we also have
  $\FSexp(\CC)=\ord(\theta)$ and so the result follows.
\end{proof}

\section{Cauchy Theorem for Quasi-Hopf algebras}\label{sec:Cauchy}
In \cite{EG99},  Etingof and Gelaki asked the following two
questions for a complex semisimple Hopf algebra $H$:
\begin{enumerate}
  \item If $p$ is a prime divisor of $\dim(H)$, does $p$ divides
  $\exp(H)$ ?
  \item If $\exp(H)$ is a power of a prime $p$, is $\dim(H)$ a power
  of $p$ ?
\end{enumerate}
The questions have been recently answered by Kashina,
Sommerh\"{a}user and Zhu \cite{KSZ}. They proved that $\exp(H)$ and
$\dim(H)$  have the same prime factors.\sprung

In this section, we generalize their result and prove an analog of
{\em Cauchy's Theorem} (Theorem \ref{t:cauchy}) for a complex
semisimple quasi-Hopf algebra $H$: ${\exp(H)}$ and $\dim(H)$ have
the same prime factors. If $H$ admits a simple self-dual module,
then $\dim(H)$ is even. Moreover, if $\dim(H)$ is odd, $\FSexp(H)$
and $\exp(H)$ are the same. \sprung

Let $\CC$ be strict spherical fusion category over $\BC$. Recall
that the $(n,r)$-th Frobenius-Schur indicator of any object in $V
\in \CC$ is a cyclotomic integer in $\BQ_n$. Let $m$ be the order of
the twist $\theta$ of $Z(\CC)$ associated with the pivotal structure
of $\CC$ and let $\zeta_m \in \BC$ be a primitive $m$-th root of
unity. Let $\BF$ be the smallest extension over $\BQ$ containing
$d_i$ for all $i \in \Gamma$. Note that $\Gal(\BF(\zeta_m)/\BF)$ is
isomorphic to a subgroup of $U(\BZ_m)$. We will call $\BF$ the {\em
dimension field} of $\CC$. Let $r$ be an integer
 such that $\sigma: \zeta_m \mapsto
\zeta_m^r$ defines an automorphism of $\BF(\zeta_m)/\BF$; in
particular, $r$ is relatively prime to $m$. By Theorem
\ref{t:nu_tr_formula}, $\nu_n(V) \in \BF(\zeta_m)$ for $V \in \CC$
and a positive integer $n$. Moreover,
$$
\sigma(\nu_n(V)) =\sigma\left(\frac{1}{\dim \CC}
\ptr(\theta_{K(V)}^{n})\right)=\frac{1}{\dim \CC}
\ptr(\theta_{K(V)}^{nr})=\nu_{nr}(V)\,.
$$
By \eqref{eq:nu1}, we have
\begin{equation}\label{eq:delta}
\delta_{i,0}= \sigma\left(\frac{1}{\dim \CC} \sum_{(X,e_X) \in
\hat{\Gamma}} N_i^X d(X) \w_{(X,e_X)}\right) =\frac{1}{\dim \CC}
\sum_{(X,e_X) \in \hat{\Gamma}} N_i^X d(X)
\w_{(X,e_X)}^r=\nu_r(X_i)\,.
\end{equation}
These equalities imply some congruences when $r$ is a prime.
\begin{prop}\label{p:divisiblity}
Let $\CC$ be a spherical fusion category over $\BC$ and
$m=\FSexp(\CC)$. For any prime $p$ such that $\zeta_m\mapsto
\zeta_m^p$ defines an automorphism of $\BF(\zeta_m)/\BF$, so in
particular $p\nmid m$, we have
$$
N_0^V\equiv N_0^{V^{\o p}} \mod{p}
$$
for all $V \in \CC$ where $N_i^X$ denotes the integer $\dim \CC(X_i,
X)$.
\end{prop}
\begin{proof}
Note that $V=\bigoplus_{i \in \Gamma} N_i^V X_i${, where
$N_i^V=\dim\CC(X_i, V)$}. By \eqref{eq:delta},
$$
\nu_p(V) =\sum_{i \in \Gamma} N_i^V \nu_p(X_i)= N_0^V\,.
$$
By the definition of Frobenius-Schur indicators,
  $
  \Tr(E_V^{(p)})=\nu_p(V)
  $
  is an integer. Since $E_V^{(p)}$ is an $\BC$-linear automorphism on $\CC(I, V^{\o p})$ of order
  1 or $p$, by a linear algebra argument in \cite{KSZ},
  $$
  \Tr(E_V^{(p)}) \equiv \dim \Hom_H(I, V^{\o p}) \mod{p}\,.
  $$
  Hence, the result follows.
\end{proof}
As an immediate consequence of Proposition \ref{p:divisiblity}, we
prove the following corollary which is a generalization of a result
in \cite{KSZ02}.\sprung
\begin{cor}
  Let $\CC$ be a spherical fusion category over $\BC$, $\BF$ the dimension field of $\CC$ and $m$ the Frobenius-Schur
  exponent of $\CC$.  Suppose that there exists a simple object $V \not\cong I$ of $\CC$ such that
  $V^*\cong V$.
  \begin{enumerate}
    \item[\rm (i)] If $\BF \cap \BQ_m =\BQ$, then $m$ is even and $(\dim\CC)^5/2$ is an algebraic integer.
 \item[\rm (ii)] If $\BF=\BQ$, then $\dim \CC$ is divisible by $2$.
  \end{enumerate}
\end{cor}
\begin{proof}
(i) Suppose that $2 \nmid m$. Then $\zeta_m \mapsto \zeta_m^2$
defines an automorphism of $\BF(\zeta_m)/\BF$
    as the condition $\BF \cap \BQ(\zeta_m) =\BQ$
    implies that $\Gal(\BF(\zeta_m)/\BF)=U(\BZ_m)$. It follows from Proposition
    \ref{p:divisiblity}  that
    $$
     0=N_0^V \equiv N_0^{V^{\o 2}} \mod{2}\,.
    $$
    This contradicts that $N_0^{V^{\o 2}}=\dim \CC(V, V^*)=1$. Therefore, $2 \mid m$. By
    \cite{Etingof02}, $\dim Z(\CC)^{5/2}/m$ is an algebraic integer. Since $\dim Z(\CC)= \dim(\CC)^2$,
    we obtain that $(\dim \CC)^5/2$ must be an algebraic integer. \sprung
(ii) If $\BF=\BQ$, then $\dim \CC$ is indeed an integer and the
assumption in (i) holds obviously. It follows from (i) that $\dim
(\CC)^5/2$ is an algebraic integer in $\BQ$ and hence an integer.
Therefore, $\dim(\CC)$ is even.
\end{proof}
Suppose further that $r$ is relatively prime to $l=\lcm(m,n)$ and
that $\sigma: \zeta_l\mapsto \zeta_l^r$ defines an automorphism of
$\BF(\zeta_l)/\BF$, where $\zeta_l$ is a primitive $l$-th root of 1.
Then $\sigma(\nu_n(V))=\nu_{n,r}(V)$ and
$\sigma(\nu_n(V))=\nu_{nr}(V)$. Thus, we have the following
proposition which generalizes the corresponding result in
\cite{KSZ}.\sprung
\begin{prop}
Let $n$ be any positive integer,  $l=\lcm(m,n)$ and  $r$ an integer
relatively prime to $l$.
 If the assignment
$\sigma:\zeta_l\mapsto \zeta_l^r$ defines an automorphism of
$\BF(\zeta_l)/\BF$, where $\zeta_l$ is a primitive $l$-th root of
unity, then
\begin{equation}
 \nu_{n,r}(V)= \nu_{nr}(V)
\end{equation}
for all $V \in \CC$, where $\zeta_l$ is a primitive $l$-th root of
unity.
\end{prop}
\begin{proof}
  The claim follows from the fact that $\sigma|_{\BF(\zeta_n)}(x)=\sigma|_{\BF(\zeta_m)}(x)$
  for all $x \in \BF(\zeta_n) \cap \BF(\zeta_m)$.
\end{proof}
Now,  we turn to the proof of the main theorem in this section. \sprung

\begin{thm}\nmlabel{Theorem}{t:cauchy}
  Let $\CC$ be a spherical fusion category over $\BC$ such that the Frobenius-Perron
  dimension of any simple object $X$ is an integer. Then $\FSexp(\CC)$, $\exp(\CC)$ and
  $\dim(\CC)$ have the same set of prime factors. Equivalently, if $H$ is a
  semisimple quasi-Hopf algebra over $\BC$, then $\FSexp(H)$, $\exp(H)$ and $\dim(H)$ have the same
  set of prime factors.
\end{thm}
\begin{proof}
  The equivalence of the two statements follows directly from the
  the independence of $\dim(\CC)$ and $\FSexp(\CC)$  on the choice of a
  pivotal structure on $\CC$ (cf. Remarks \ref{r:piv_dim} and \ref{r:fsexp}),
  and the characterization of $\C{H}$ for some semisimple quasi-Hopf algebra
  over $\BC$ as a fusion category over $\BC$ with integer
  Frobenius-Perron dimension for each simple object by \cite[section
  8]{ENO}. \sprung

  Without loss of generality, we may assume $\CC$ is strict
  pivotal. Since both the Frobenius-Schur exponent and the
  dimension of $\CC$ are independent of the choice of a spherical
  structure, we can assume that $\CC$ is endowed with the
  canonical spherical structure with $d(V)$ equal to the
  Frobenius-Perron dimension of $V \in \CC$. The center $Z(\CC)$ is a modular tensor category
  with the twist $\theta$ associated with the canonical spherical structure of $\CC$.
  By \cite{Etingof02}, the order of $\theta$ divides
  $\dim(Z(\CC))^{5/2}=\pm \dim(\CC)^5$. Thus, the prime factors of
  $\ord(\theta)$ are also prime factors of $\dim(\CC)$.
  Let $\{X_i\}_{i \in \Gamma}$ be a complete set of non-isomorphic  simple objects of
  $\CC$. Recall from \cite{Etingof02} that  $V=\sum_{i \in \Gamma}d(X_i) X_i$ defines
  a rank one ideal in the Grothendieck ring of $\CC$ such that $X_i \o V \cong d(X_i)V$ and
  $$
  d(V) =  \sum_{i \in \Gamma} d(X_i)^2
  =\dim (\CC)\,.
  $$
  Suppose $p$ is not a prime factor of
  $\ord(\theta)$. Since $d(X_i)$ is an integer for all $i \in \Gamma$,
  $\BF=\BQ$.
  By Proposition \ref{p:divisiblity}, we have
  $$
   N_0^V \equiv N_0^{V^{\o p}} \mod{p}\,.
  $$
  Since
  $$
  V^{\o p} \cong d(V^{\o (p-1)}) V = d(V)^{p-1} V  = \dim(\CC)^{p-1}
  V,
  $$
  we find
  $$
  \dim \CC(I, V^{\o p})  = \dim(\CC)^{p-1}\,.
  $$
  Therefore,
  $$
  1 \equiv \dim(\CC)^{p-1} \mod{p}
  $$
  and hence $p \nmid \dim(\CC)$. Thus, by Theorem
  \ref{t:FS_exponent},   $\FSexp(\CC)$ and $\dim(\CC)$ have the same set of prime
  factors. \sprung

  By Corollary \ref{c:6.2},
  $$
  \FSexp(\CC)=\exp(\CC)\quad\text {or}\quad \FSexp(\CC)=2\exp(\CC).
  $$
  To complete the proof, it suffices to show that $\exp(\CC)$ is even whenever $\FSexp(\CC)$ is
  even. Suppose $\FSexp(\CC)$ is even and $\exp(\CC)$ is odd. Then we have
  $\FSexp(\CC)=2 \exp(\CC)$. Recall that $\FSexp(\CC)$ is the order
  of $\theta$. By \nmref{l:change_twist}, there exists another spherical
  pivotal structure $\hat{j}$ on $Z(\CC)$ such that
  $\ord(\hat{\theta})=\exp(\CC)$ for the associated twist $\hat\theta$, while
  $\hat{d}(V)=\pm d(V)$
  for any simple object $V$ of $Z(\CC)$, where $\hat{d}$ is the dimension
  function associated with $\hat{j}$. Let $\widehat{\CC}$ denote
  the spherical fusion category $Z(\CC)$ endowed with the spherical structure
  $\hat{j}$. By \nmref{l:change_twist}, $\hat{\CC}$ is modular since $Z(\CC)$ is a MTC. By
  Theorem \ref{t:o(theta)}, Corollary \ref{c:invariance}, and Remark \ref{r:fsexp}, we have
  $$
   \FSexp(\CC) = \FSexp(Z(\CC))=\FSexp(\hat{\CC}) =
   \ord(\hat{\theta}) =\exp(\CC)\,,
  $$
  a contradiction!
\end{proof}

\begin{cor}
Let $H$ be an odd dimensional semisimple quasi-Hopf algebra over
$\BC$. Then the exponent of $H$ and Frobenius-Schur exponent of $H$
are identical.
\end{cor}
\begin{proof}
 By the preceding theorem, $\FSexp(H)$ is odd and so the
claim follows from Corollary \ref{c:6.2}.
\end{proof}

\section{Bounds on the Exponent}\label{sec:bounds}

A direct application of Etingof's bound on the order of the twist of
a modular tensor category \cite{Etingof02} shows that the
Frobenius-Schur exponent (and hence the exponent) of a spherical
fusion category divides the fifth power of its dimension. In this
section we will strengthen this bound in important special cases.
For a semisimple quasi-Hopf algebra $H$, we can show that the
Frobenius-Schur exponent of $H$ divides $\dim(H)^4$. The techniques
in this case are close to those of Etingof and Gelaki in
\cite{EG02}. Special care has to be taken to deal with the
nontrivial associativity isomorphisms (which can be avoided in
\cite{Etingof02} by using categorically defined determinants, at the
cost of a higher bound). Our bound is higher than the bound of
$\dim(H)^3$ obtained for Hopf algebras by Etingof and Gelaki, which
can perhaps be tracked to the fact that we cannot use the dual Hopf
algebra $H^*$. For the special but important class of
group-theoretical quasi-Hopf algebras introduced by Ostrik
\cite{Ostrik03}, on the other hand, we derive a bound of
$\dim(H)^2$. This bound, which is even better than the best
previously known bound for general semisimple \emph{ordinary} Hopf
algebras, was recently obtained by Natale \cite{Nat05}. Using the
general theory of the Frobenius-Schur exponent, we can reduce the
problem to the case where the quasi-Hopf algebra is just a dual
group algebra with a quasi-bialgebra structure induced by a
three-cocycle. In this case, we can compute the Frobenius-Schur
exponent directly by considering the indicators of individual
representations, or really (after a further reduction using the
invariance properties) of just one example.\sprung

\begin{thm}
  Let $H$ be a semisimple complex quasi-Hopf algebra. Then the
  Frobenius-Schur exponent of $H$ divides $\dim(H)^4$.
\end{thm}
\begin{proof}
  For $X\in\C H$ and $V\in\C{D(H)}= Z(\C H)$ we write
  $$\hat e_V(X)=\tau_{XV}e_V(X)\colon V\o X\to V\o X,$$ where $\tau$
  is the ordinary vector space flip. From the hexagon equation
  $$e_V(X\o Y)=\Phi_{XYV}\inv(X\o e_V(Y))\Phi_{XVY}(e_V(X)\o
  Y)\Phi_{VXY}\inv$$
  we deduce
  $$\det\hat e_V(X\o Y)=\det (\Phi\inv_{XYV})\det(\hat e_V(Y))^{\dim
  X}\det(\Phi_{XVY})\det(\hat e_V(X))^{\dim(Y)}\det(\Phi\inv_{VXY})$$
  Specializing $Y=H$ and using $X\o H\cong H^{\dim X}$, we find
  \begin{multline*}
    \det(\hat e_V(H))^{\dim X}=\det(\hat e_V(X\o H))
    \\=\det(\Phi\inv_{XHV})\det(\hat e_V(H))^{\dim
    X}\det(\Phi_{XVH})\det(\hat e_V(X))^{\dim H}\det(\Phi\inv_{VXH})
  \end{multline*}
  and thus
  \[\det(\hat e_V(X))^{\dim H}=\det(\Phi_{XHV})\det(\Phi_{VXH})\det(\Phi_{XVH}\inv).\]
  For $V,W\in\C{D(H)}$ we then find
  \begin{multline*}
   \det(e_V(W)e_W(V))^{\dim H}
     =\det(\hat e_V(W))^{\dim H}\det(\hat e_W(V))^{\dim H}
\\=\det(\Phi_{WHV})\det(\Phi_{VWH})\det(\Phi\inv_{WVH})\det(\Phi_{VHW})\det(\Phi_{WVH})\det(\Phi\inv_{VWH})
     \\=\det(\Phi_{WHV})\det(\Phi_{VHW}).
  \end{multline*}
  The other hexagon equation
  \[e_{V\o W}(X)=\Phi_{XVW}(e_V(X)\o W)\Phi\inv_{VXW}(V\o e_W(X))\Phi_{VWX}\]
  gives
  \[\det(\hat e_{V\o W}(X))=\det(\hat e_V(X))^{\dim W}\det(\hat e_W(X))^{\dim V}
    \det(\Phi_{VWX})\det(\Phi\inv_{VXW})\det(\Phi_{XVW}).\]
  For $W=D(H)$, using $V\o D(H)\cong D(H)^{\dim V}$ as
  $D(H)$-modules, and abbreviating $D:=D(H)$, we find
  \begin{align*}
    \det(\hat e_{D}(X))^{\dim V}
      &=\det(\hat e_V(X))^{\dim D}\det(\hat e_{D}(X))^{\dim
      V}\det(\Phi_{VDX})\det(\Phi\inv_{VXD})\det(\Phi_{XVD})
  \end{align*}
  and hence
  $$
  \begin{aligned}
    \det(\hat e_V(X))^{-\dim D}
      &=\det(\Phi_{VDX})\det(\Phi\inv_{VXD})\det(\Phi_{XVD})\\
      &=\det(\Phi_{VHX})^{\dim H}\det(\Phi\inv_{VXH})^{\dim
      H}\det(\Phi_{XVH})^{\dim H}
  \end{aligned}
  $$
  because the associativity isomorphisms depend only on the
  $H$-module structures of the objects involved, and $D(H)\cong
  H^{\dim H}$ as $H$-module. Comparing with the previous
  calculation, which gives
  \[\det(\hat e_V(X))^{\dim D}
    =\left(\det(\hat e_V(X))^{\dim H}\right)^{\dim H}
    =\det(\Phi_{XHV})^{\dim H}\det(\Phi_{VXH})^{\dim H}\det(\Phi_{XVH}\inv)^{\dim H}\]
  we find
  \[\det(\Phi_{VHX})^{\dim H}\det(\Phi_{XHV})^{\dim H}=1,\]
  and in particular
  \[\det(e_V(D)e_{D}(V))^{\dim H}=\det(\Phi_{DHV})\det(\Phi_{VHD})
    =\det(\Phi_{HHV})^{\dim H}\det(\Phi_{VHH})^{\dim H}=1.\]
  We continue as in \cite{Etingof02} and \cite{EG02}:
  Since $\theta_{V \o D(H)}$ and $\theta_V\o\theta_{D(H)}$ agree up
  to a factor $e_V(D)e_D(V)$, we have
  \begin{multline*}\det(\theta_{D(H)})^{\dim V\dim H}=\det(\theta_{V \o D(H)})^{\dim
  H}
  =\det(\theta_V\o\theta_{D(H)})^{\dim H} \\=
  \det(\theta_V)^{\dim D(H)\dim H}\det(\theta_{D(H)})^{\dim V\dim H}
  \end{multline*}
  and thus $\det(\theta_V)^{\dim D(H)\dim H}=1$. If $V$ is simple, and thus
  $\theta_V$ is a scalar, we conclude
  $$
  \theta_V^{\dim V\dim H^3}=1.
  $$
  Since $\dim V$ divides $\dim H$,
  we are done.
\end{proof}

Now we turn to the announced bound on the (Frobenius-Schur) exponent
of a group-theoretical quasi-Hopf algebra.\sprung

Let $G$ be a finite group and $\omega\colon G^3\to \BC^\times$ a
three-cocycle. Let $H\subset G$ a subgroup, and $\psi$ a 2-cochain
of $G$ with $d\psi=\omega|_H$. The group-theoretical category
$\CC(G,H,\omega,\psi)$ defined by Ostrik \cite{Ostrik03} is the
category of $\BC_{\psi}[H]$-bimodules in the category
$\CC(G,\omega)$, which in turn is the category of $G$-graded vector
spaces, made into a monoidal category with the usual tensor product
but nontrivial associativity constraint $\Phi$ given by $\omega$.
More precisely $\Phi_{UVW}\colon (U\o V)\o W\to U\o(V\o W)$ is given
by $\Phi(u\o v\o w)=\omega(|u|,|v|,|w|)u\o v\o w$ if $u,v,w$ are
homogeneous (and $|x|$ denotes the degree of a homogeneous element
$x$). It makes sense to consider $\BC_{\psi}[H]$-bimodules in
$\CC(G,\omega)$ since the twisted group algebra $\BC_{\psi}[H]$ is
an associative algebra in the monoidal category $\CC(G,\omega)$
(thanks to the condition $d\psi=\omega|_H$). By
\cite{Schauenburg01}, the center of $\CC(G,H,\omega,\psi)$ is
isomorphic, as a braided monoidal category, to the center of
$\CC(G,\omega)$. The quasi-Hopf algebras over $\BC$, whose
representation categories are monoidally equivalent to
group-theoretical categories, are called group-theoretical
quasi-Hopf algebras. The center of the category $\CC(G,\omega)$ is
isomorphic to the category of modules over the twisted Drinfeld
double $D^\omega(G)$ of Dijkgraaf, Pasquier, and Roche \cite{DPR90}.
Thus we see that for any group-theoretical quasi-Hopf algebra $K$,
there exist a group $G$ and a 3-cocycle $\omega$ on $G$ such that
the module categories over $D(K)$ and $D^\omega(G)$ are equivalent
braided monoidal categories, that is, $D(K)$ and $D^\omega(G)$ are
gauge equivalent quasitriangular quasi-Hopf algebras. This argument
can be found in \cite{Natale03}, where Natale also shows that
group-theoretical quasi-Hopf algebras can in fact be characterized
as those whose double is gauge equivalent to a twisted double of a
finite group.\sprung

The Frobenius-Schur indicators of the objects of $\CC(G,\omega)$
were already computed in \cite{NS052}. In the terminology of
Frobenius-Schur indicators, \cite[Prop. 7.1]{NS052} says that the
Frobenius-Schur exponent of a simple object $V_x$ associated with
the element $x \in G$ is equal to $\ord(x)\ord(\res_{\langle
x\rangle}[\w])$. Hence $\FSexp(\CC(G,\omega))$ is the least common
multiple of the numbers $\ord(x)\ord(\res_{\langle x\rangle}[\w])$
for all $x\in G$ (see also below). But by Corollary
\ref{c:invariance} we know that the Frobenius-Schur exponent is
invariant under the center construction. Hence, we have proved the
following result: \sprung
\begin{thm}
  The Frobenius-Schur exponent of the group-theoretical category
  $\CC(G,H,\omega,\psi)$ is the least common multiple of the numbers
  $|C|\cdot \ord(\res_C[\w])$ where $C$ runs through the (maximal)
  cyclic subgroups of $G$. In particular, $\FSexp(\CC(G,H,\omega,\psi))$
  divides $\exp(G)^2$ and
  $\exp(G) \ord([\omega])$. Hence, for each group-theoretical quasi-Hopf algebra $K$,
  we have $\exp(K)|\FSexp(K)|\dim(K)^2$. \sprung
\end{thm}

The result on the Frobenius-Schur exponent of $\CC(G,\omega)$ cited
above relies on the fact that $\CC(G,\omega)$ can be described as
the category of modules over a certain quasi-Hopf algebra
$H(G,\omega)$. The proof of \cite[Prop.7.1]{NS052} specializes
general formulas for indicators over quasi-Hopf algebras to this
case. It may be interesting to see a proof that computes the
indicators from scratch using their first definition in \cite{NS05}
and the description of $\CC(G,\omega)$ above. We will do this in the
rest of the section.\sprung

A complete set of representatives for the isomorphism classes of
simple objects in $\CC(G,\omega)$ is given by $V_g={\BC}$ as a
vector space, made into a homogeneous graded vector space of degree
$g$, for each $g\in G$. We will treat the canonical isomorphisms
$V_g\o V_h\cong V_{gh}$ as identities. As a consequence, every
morphism between iterated tensor products of simples is given as
multiplication with a scalar, and we will sometimes identify the
morphism and the scalar below. It remains to make a suitable choice
of dual objects. We will take $(V_g)\du:=(V_{g\inv},\ev_g,\db_g)$
with $\db_g=1:{\BC}\to V_g\o V_{g\inv}$, which forces us to choose
$\ev_g=\omega(g,g\inv,g)\inv=\omega(g\inv,g,g\inv)$, since we have
to ensure that
\[V_g\xrightarrow{\db_g\o V_g}(V_g\o V_{g\inv})\o V_g
\xrightarrow{\Phi}V_g\o(V_{g\inv}\o
V_g)\xrightarrow{V_g\o\ev_g}V_g\] is the identity, and
$\Phi_{V_g,V_{g\inv},V_g}=\omega(g,g\inv,g)$.\sprung

We proceed to determine the pivotal structure of $\CC$. The
component $j_g\colon V_g\to (V_g)\bidu=V_g$ is determined by the
requirement that the composition
\[{\BC}\xrightarrow{\db_g}V_g\o V_{g\inv}\xrightarrow{j_g\o V_{g\inv}}V_g\o V_{g\inv}\xrightarrow{\ev_{g\inv}}{\BC}\]
be the identity. Thus $j_g=\omega(g\inv,g,g\inv)$, since
$\ev_{g\inv}=\omega(g,g\inv,g)$.\sprung

Finally, let us investigate the explicit form of the isomorphism
\[D\colon \CC(I,V\o W)\to\CC(I,W\o V\bidu)\]
introduced in \cite[Definition 3.3]{NS05}. We have
\[D(f)=\left(\BC\xrightarrow{\db_{V\du}}V\du\o V\bidu
  \xrightarrow{V\du\o f\o V\bidu}(V\du\o (V\o W))\o V\bidu
  \xrightarrow{\Theta\o V\bidu}W\o V\bidu\right),\]
where $\Theta=\left((V\du\o (V\o W))
  \xrightarrow{\phi\inv}((V\du\o V)\o W)
  \xrightarrow{\ev_V\o W}W\right).$ Note that if $W=V\du$, then
\[\left(V\du\xrightarrow{V\du\o\db_V}V\du\o (V\o V\du)\xrightarrow{\Theta}V\du\right)=\id_{V\du}.\]
Now if $V=V_g$ and $W=V_{g\inv}$, then all the morphisms involved
can be identified with scalars. We have $\db_{V_g}=1$, and so
$\Theta=1$, and since also $\db_{V_{g\inv}}=1$, we finally have
$D(f)=f$. In particular, we see that the isomorphism
$E_{V_g,V_{g\inv}}\colon\CC({\BC},V_g\o{V_{g\inv}})\to\CC({\BC},V_{g\inv}\o
V_g)$ is given, under the identification of both sides with
$\CC({\BC},{\BC})\cong {\BC}$, by the scalar
$j_g\inv=\omega(g,g\inv,g)$. \sprung

Let us now calculate the higher indicators of $V_g$ in the special
case $G=\Z/N\Z$ and $g=\ol 1$. To describe, as in \cite{MS89}, the
cocycles on $G$, we define $\hh n\in\{0,\dots,N-1\}$ for $n\in\Z$ by
$\hh n\equiv n$ modulo $N$. Then the class of the cocycle $\omega_1$
defined by
\[\omega_1(\ol\ell,\ol m,\ol n)=\exp\left(\frac{2\pi i}{N^2}\hh\ell(\hh m+\hh n-\hh{m+n})\right)\]
generates the group $H^3(\Z/N\Z,{\BC^\times})\cong C_N$. In
particular, every cocycle on $\Z/N\Z$ has the form
$\omega_t=\omega_1^t$ for some $0\leq t<N$, with
\[\omega_t(\ol\ell,\ol m,\ol n)=\exp\left(\frac{2\pi it}{N^2}\hh\ell(\hh m+\hh n-\hh{m+n})\right).\]
We will determine the indicators of $V_{\ol
1}\in\CC(\Z/N\Z,\omega_t)$. We start by noting that
\[\omega_t(\ol1,\ol n,\ol1)=\exp\left(\frac{2\pi it}{N^2}(\hh n+1-\hh{n+1})\right)
 =\begin{cases}
 \exp\left(\frac{2\pi it}{N}\right)=\zeta^t&\text{ if }\ol n=\ol{N-1}\\
 1&\text{ otherwise,}
 \end{cases}\]
where we have put $\zeta=\exp\left(\frac{2\pi i}{N}\right)$. In
particular, $E_{V_{\ol 1},V_{\ol{N-1}}}$ is multiplication with
$\zeta^t$. We can also compute the scalar $\Phi_n\colon V_{\ol
1}^{\o (n-1)}\o V_{\ol 1}\to V_{\ol 1}^{\o n}$. We have the
recursion formula $$\Phi_n=(V_{\ol 1}\o\Phi_{n-1})\Phi_{V_{\ol
1},V_{\ol 1}^{\o(n-2)},V_{\ol 1}}=\Phi_{n-1}\omega_t(\ol
1,\ol{n-2},\ol 1),$$ and hence $\Phi_{\ell N}=\zeta^{t(\ell-1)}$.
Finally, by definition of $E_{V_{\ol 1}}^{(n)}$, which is again a
scalar or zero, we know that $\nu_n(V_{\ol 1})=0$ if $N\nmid n$,
while $\nu_{\ell N}(V_{\ol
1})=\zeta^{t}\phi_{N\ell}=\zeta^{t\ell}$.\sprung

As a consequence of these calculations, we find:\sprung
\begin{lem}
Let $G$ be a finite group, and $\omega\colon G^3\to{\BC^\times}$ a
three-cocycle representing $[\omega]\in H^3(G,{\BC^\times})$.
\begin{enumerate}
\item
  If $G$ is a cyclic group of order $N$, then
  $\FSexp(\CC(G,\omega))=N\cdot\ord([\omega])$. In particular we have
  $N|\FSexp(\CC(G,\omega))|N^2$, and for each
  $n|N$ there is $\omega$ with $\FSexp(\CC(G,\omega))=Nn$.
\item ${\FSexp(\CC(G,H,\omega,\psi))}$ is the least common
multiple of the
  numbers
  $\ord(C)\ord(\res_C[\omega])$ where $C$ runs through the (maximal)
  cyclic subgroups of $G$.
\end{enumerate}
\end{lem}
\begin{proof}
  The calculations preceding the statement of the lemma have shown
  that $\FSexp(V_{\ol 1})=N\ord(\zeta^t)=N\ord([\omega])$ if $G=\Z/N\Z$ and
  $\omega=\omega_t$. But since the left hand side also only depends on the cohomology
  class of $\omega$, the formula holds for all cocycles. Thus, more generally,
  $\FSexp(V_g)=N\ord([\omega])$ whenever $g$ generates the cyclic group
  $G$, because the group isomorphism $f\colon G\to\Z/N\Z$ mapping $g$ to
  $\ol 1$ induces a monoidal category equivalence
  $\CC(G,\omega)\to\CC(\Z/N\Z,\omega')$ mapping $V_g$ to $V_{\ol 1}$,
  where $\omega'$ corresponds to $\omega$ under the isomorphism
  $H^3(\Z/N\Z,{\BC^\times})\to H^3(G,{\BC^\times})$ induced by $f$.
  For general $g\in G$, we may compute the Frobenius-Schur indicators and exponent
  of $V_g\in\CC(G,\omega)$ by restricting ourselves to the full monoidal subcategory
  $\CC(\langle g\rangle,\omega|_{\langle g\rangle^3})$ containing $V_g$. Thus
  $$\FSexp(V_g)\mid \ord(g)\ord([\omega_{\langle g\rangle^3}])\mid N\ord([\omega]).$$
  Summing up, the Frobenius-Schur exponents of all the simples
  $V_g$ divide $N\ord([\omega])$, and equality occurs for any $g$ generating
  $G$. In particular $\FSexp(\CC(G,\omega))=N\ord([\omega])$, and
  of course there is a cohomology class of order $n$ for any
  divisor of $N$ since $H^3(G,{\BC^\times})\cong C_N$.\sprung

For a general finite group $G$ and  any simple
$V_g\in\CC(G,\omega)$, we can compute its Frobenius-Schur exponent
inside the category $\CC(\langle g\rangle,\omega|_{\langle
g\rangle^3})$, and in particular inside one of the categories
$\CC(C,\omega|_{C^3})$ for $C$ a cyclic subgroup of $G$, which we
may choose to be a maximal cyclic subgroup. The Frobenius-Schur
exponent is thus the least common multiple of the Frobenius-Schur
exponents of these categories, which we have already computed.
\end{proof}

\bibliographystyle{amsalpha}
 \providecommand{\bysame}{\leavevmode\hbox
to3em{\hrulefill}\thinspace}
\providecommand{\MR}{\relax\ifhmode\unskip\space\fi MR }
\providecommand{\MRhref}[2]{%
  \href{http://www.ams.org/mathscinet-getitem?mr=#1}{#2}
} \providecommand{\href}[2]{#2}

\end{document}